\documentclass[12pt]{article}

\usepackage{latexsym,amsmath,amscd,amssymb,graphics}
\usepackage{enumerate}

\usepackage{graphicx}

\usepackage[square,authoryear]{natbib}
\usepackage[colorlinks]{hyperref}
\usepackage{url}

\usepackage[all]{xy}

\makeatletter

\@addtoreset{figure}{section}
\def\thefigure{\thesection.\@arabic\c@figure}
\def\fps@figure{h, t}
\@addtoreset{table}{bsection}
\def\thetable{\thesection.\@arabic\c@table}
\def\fps@table{h, t}
\@addtoreset{equation}{section}

\makeatother

\begin{document}

\newtheorem{theorem}{Theorem}[section]
\newtheorem{definition}[theorem]{Definition}
\newtheorem{lemma}[theorem]{Lemma}
\newtheorem{remark}[theorem]{Remark}
\newtheorem{proposition}[theorem]{Proposition}
\newtheorem{corollary}[theorem]{Corollary}
\newtheorem{example}[theorem]{Example}

\def\below#1#2{\mathrel{\mathop{#1}\limits_{#2}}}



\title{Group Actions on Chains of
Banach Manifolds and Applications to Fluid Dynamics}
\author{Fran\c{c}ois Gay-Balmaz$^{1}$ and Tudor S. Ratiu$^{1}$}
\addtocounter{footnote}{1} \footnotetext{Section de
Math\'ematiques, \'Ecole Polytechnique F\'ed\'erale de Lausanne.
CH--1015 Lausanne. Switzerland.
\texttt{Francois.Gay-Balmaz@epfl.ch, Tudor.Ratiu@epfl.ch}
\addtocounter{footnote}{1} }

\date{ }
\maketitle

\makeatother
\maketitle




\begin{abstract}
This paper presents the theory of non-smooth Lie group actions on
chains of Banach manifolds. The rigorous functional analytic
spaces are given to deal with quotients of such actions. A
hydrodynamical example is  studied in detail.
\end{abstract}



\section{Introduction}
\label{section: Introduction} The goal of this paper is to show
the existence of \textit{slices}, to study the geometric
properties of the \textit{orbit type sets}, the \textit{fixed
point sets}, the \textit{isotropy type sets}, and the \textit{orbit
spaces} associated to a certain important class of non-smooth
actions of a Lie group $G$ on a chain of Banach manifolds. These
considerations will be applied to the motion of an ideal
homogeneous incompressible fluid in a compact domain.

This study is motivated by many examples that appear as results of
various reduction procedures of infinite dimensional Hamiltonian
systems with symmetry. All these applications have two technical
difficulties in common. First, the evolutionary equations do not
have smooth flows in the function spaces that are natural to the
problem. If the system is linear, this corresponds to the fact
that the equation is defined by an unbounded operator. This leads
to the study of vector fields on manifolds of maps that are
defined only on open dense subsets (see \cite{ChMa1974} for a
presentation of this theory). Second, the actions, even those of
finite dimensional Lie groups, are not smooth in the usual sense.
The dependence on the group variable is only smooth on an open
dense set. It turns out that in spite of these problems one can
endow various objects related to the symmetry of the manifold and
the flows of various interesting vector fields with certain weak
smooth structures that are compatible between themselves. It is
the goal of this paper to present the beginnings of such a theory
and work out a concrete example coming from fluid dynamics.

Before presenting the outline of this theory and discuss some
applications let us recall briefly some standard results about the
structure of orbit spaces in the finite dimensional case.
Consider a smooth and proper action
\[
\Phi : G\times M\longrightarrow M
\]
of a finite dimensional Lie group $G$ on a smooth finite
dimensional manifold $M$. Then we have the following results.
\begin{enumerate}
\item[{\rm (i)}] If the action is free, $M/G$ is a smooth manifold
and the canonical projection $\pi : M\longrightarrow M/G$ defines
on $M$ the structure of a smooth principal $G$-bundle. \item[{\rm
(ii)}] If all the isotropy subgroups are conjugate to a given one,
say $H\subset G$, then $M/G$ is a smooth manifold and the
canonical projection $\pi : M\longrightarrow M/G$ defines on $M$
the structure of a smooth locally trivial fiber bundle with
structure group $N(H)/H$ and fiber $G/H$, where
\[
N(H):=\{g\in G\,|\,gHg^{-1}=H\}
\]
is the normalizer of $H$ in $G$. \item[{\rm (iii)}] For an
arbitrary smooth proper action, if $H\subset G$ is a closed subgroup of
$G$, define the following subsets of $M$:
\begin{align*}
M_{(H)}&:=\left\{m\in M\,|\,G_m\text{ is conjugate to } H \right\}, \\
M^H&:=\left\{m\in M\,|\, H \subset G_m\right\}, \\
M_{H}&:=\left\{m\in M\,|\,G_m=H\right\}.
\end{align*}
$M_{(H)}$ is called the $(H)$-\textbf{orbit type set}, $M^H$ is
the $H$-\textbf{fixed point set}, and $M_H$ is the
$H$-\textbf{isotropy type set}. $M_{(H)}$ is $G $-invariant whereas
$M_H$ and $M^H$ are not $G$-invariant, in general. In addition, $M_H
\subset M^H$.

\item[{\rm (iv)}] $M_{(H)}, M^H$, and $M_H$ are submanifolds of $M$.
Moreover,
$M_H$ is open in $M^H$ and for $m\in M^H$,
\[
T_mM^H=\left\{v_m\in T_mM\,|\,T_m\Phi_h(v_m)=v_m,\,\forall\,h\in
H\right\}.
\]
We have the partitions
\[
M=\bigsqcup_{(H)} M_{(H)}\,\text{ and
}\,M/G=\bigsqcup_{(H)}(M_{(H)}/G);
\]
$M_{(H)}/G$ is called the \textbf{isotropy stratum of type} $(H)$.
By (ii), $\pi_{(H)}:=\pi|_{(H)}:M_{(H)}\longrightarrow M_{(H)}/G$
is a smooth locally trivial fiber bundle with structure group
$N(H)/H$ and fiber $G/H$.
\end{enumerate}
Of course (i) is a particular case of (ii), and (ii) is a
particular case of (iii). With the additional hypothesis that the
orbit map $\Phi^m :G \rightarrow M $ is an immersion for all $m
\in M $, the result in (i) is still valid for Banach Lie groups
acting smoothly and properly on Banach manifolds (see
\cite{Bo1972}, Ch III, \S 1, Proposition 10) . This hypothesis
always holds in finite dimensions.

We now present some examples of $G$-actions on Banach manifolds
that arise in the reduction by stages procedure of infinite
dimensional dynamics. We shall see that the previous results are
not applicable to these actions since they are not
smooth.\\

\noindent\textbf{First example : incompressible fluid dynamics}\\

Let $(M,g)$ be a smooth Riemannian manifold and denote by
$Iso: = Iso(M,g)$ the finite dimensional Lie group of
isometries of $(M,g)$ (the \cite{MeSt1939} Theorem). The Lie group
topology of $Iso$ coincides with the topology of uniform convergence
on compact sets. If $M $ is compact, the Lie group $Iso$, and
hence $Iso^+$, is a compact Lie group. For the proof of these statements
see \cite{KoNo1963}, Theorem 3.4 in Chapter VI.

We consider the motion of an ideal incompressible fluid in a
compact oriented Riemannian manifold $M$ with boundary. The
appropriate configuration space is $\mathcal{D}_\mu^s (M),
s>\frac{\operatorname{dim}M}{2}+1$, the Hilbert manifold of volume
preserving $H^s$-diffeomorphisms of $M$. The Lagrangian is the
quadratic form associated to the weak $L^2$ Riemannian
metric on $\mathcal{D}_\mu^s (M)$ given by
\[
\langle\!\langle u_\eta,v_\eta
\rangle\!\rangle=\int_Mg(\eta(x))(u_\eta(x),v_\eta(x))\mu
(x),\;\;u_\eta,v_\eta\in
T_\eta \mathcal{D}^s_\mu(M),
\]
where $g$ is the Riemannian metric and $\mu$ is the associated Riemannian
volume form.

Since this Lagrangian is invariant under the
following two commuting actions
\[
R : \mathcal{D}_\mu^s(M)\times T\mathcal{D}_\mu^s(M)
\longrightarrow T\mathcal{D}_\mu^s(M),\;\;R_\eta(v_\xi)=v_\xi\circ
\eta\,\,\text{ and}
\]
\[
L : Iso^+\times T\mathcal{D}_\mu^s(M)\longrightarrow T\mathcal{D}
_\mu^s(M),\;\;L_i(v_\xi)=Ti\circ v_\xi,
\]
it is formally true that the Poisson reduction by stages procedure
can be applied, that is, the reduction by $\mathcal{D}_\mu^s (M)$ and
then by $Iso^+$ coincides with the one step reduction by the
product group $\mathcal{D}_\mu^s(M) \times Iso^+$.

The reduction by $\mathcal{D}_\mu^s(M)$ (first stage reduction) is
well known and leads to the Euler equations for an ideal
incompressible fluid on the first reduced space
$\mathfrak{X}^s_{div}(M)=T\mathcal{D}_\mu^s (M)/
\mathcal{D}_\mu^s(M)$.

Our goal is to carry out in a precise sense the reduction by
$Iso^+$ (second stage reduction). In spite of the fact that
$Iso^+$ is a compact finite dimensional Lie group, several problems
occur. The action $l$ of $Iso^+$ on $\mathfrak{X}^s_{div}(M)$
induced by $L$ is given by
\[
l : Iso^+\times \mathfrak{X}^s_{div}(M)\longrightarrow
\mathfrak{X} ^s_{div}(M),\;\;l_i(u)=Ti\circ u\circ i^{-1}=i_*u
\]
and we remark that $l$ is neither free nor $C^1$ ($l$ is $C^1$ as
a map with values in $\mathfrak{X}^{s-1}_{div}(M)$). Thus the usual
results about the geometric properties of the isotropy type submanifolds
and the orbit space, valid for the case of a smooth proper
action, cannot be used here. So it would be interesting to
determine the exact differentiable structure of
$\mathfrak{X}^s_{div}(M)_H$, $\mathfrak{X}^s_{div}(M)_H/N(H)$,
$\mathfrak{X}^s_{div}(M)_{(H)}$,
$\mathfrak{X}^s_{div}(M)_{(H)}/Iso^+$ in order to define tangent
bundles, vector fields, evolution equations on them, and to
carry out the second stage reduction procedure.

The following relevant facts are a guide in the search of a
useful definition for the tangent bundle of the orbit spaces:\\
$(1)$ the infinitesimal generators of the action $l$ are not
vector fields on $\mathfrak{X}^s_{div}(M)$ but they are sections
of the vector bundle
$\mathfrak{X}^s_{div}(M)\times\mathfrak{X}^{s-1}_{div}(M)
\longrightarrow \mathfrak{X}^s_{div}(M)$,\\
$(2)$ the Hamiltonian vector field is given by
$X_h(u)=(u,-\operatorname{P}_e (\nabla_uu))$, so it is not a
vector field on $\mathfrak{X}^s_{div}(M)$, since it takes values
in $\mathfrak{X}^s_{div}(M)\times\mathfrak{X}^{s-1}_{div}
(M)$.\\
$(3)$ the integral curves of the Euler equation are not $C^1$
curves but are elements of $C^0(I,\mathfrak{X}^s_{div}(M))\cap
C^1(I,\mathfrak{X}^{s-
1}_{div}(M))$.\\
We will treat this example in detail in Section $6$. Note that
the same situation arises in
the  motion of the averaged incompressible fluid.\\

\noindent\textbf{Second example : a nonlinear wave equation}\\

On the configuration space $H^s(S^1,\mathbb{R}^2)\times
H^{s-1}(S^1,\mathbb{R} ^2), s\geq 2$, we consider the Lagrangian
\[
L(\varphi,\dot\varphi)=\frac{1}{2}\langle
\dot\varphi,\dot\varphi\rangle- \frac{1}{2}\langle
\varphi',\varphi'\rangle+\frac{1}{4}\langle\varphi,\varphi\rangle^2
\]
where $\langle\,,\rangle$ is the $L^2$ inner product, the dot denotes
time derivative, and the prime space derivative. Then the integral curves
$\varphi(t)$ of the corresponding Euler-Lagrange equation are in fact
periodic solutions of the nonlinear wave equation
\[
\ddot\varphi(t)=\varphi(t)''+\varphi(t)|\varphi(t)|^2.
\]
The lagrangian $L$ is invariant under the two following commuting
actions:
\[
L : SO(2)\times (H^s\times H^{s-1})\longrightarrow H^s\times H^{s-
1},\,\,\,L_A(\varphi,\psi)(s)=(A\!\cdot \varphi(s),A\!\cdot
\psi(s))
\]
\[
R : S^1\times (H^s\times H^{s-1})\longrightarrow H^s\times H^{s-
1},\,\,\,R_\alpha(\varphi,\psi)(s)=(\varphi(s+\alpha),\psi(s+\alpha)).
\]
As before, it should be formally possible to apply the Poisson
reduction by stages. Since the action of $SO(2)$ on
$P^s:=H^s\times H^{s-1}-\{(0,0)\}$ is smooth free and proper, the
first reduced space $P^s/SO(2)$ is a Hilbert manifold (see
\cite{Bo1972}, Chapter III, Proposition 10).

Some problems occur for the second reduced space. In fact, as
was the case in the first example, the action of $S^1$ on $P^s/SO
(2)$ induced by $R$ is neither free nor $C^1$. The same problems
arise in the symplectic reduction by stages where $S^1$ acts on
$\textbf{J}^{-1}_{SO(2)} (\mu)/SO(2)$ or on
$\textbf{J}_{\mu,S^1}^{-1}(\nu)$. Here $\textbf{J}_{SO (2)},
\textbf{J}_{S^1}$ denote the momentum mappings for the
corresponding actions, and $\textbf{J}_{\mu,S^1}$ denotes the map
induced by $\textbf
{J}_{S^1}$ on $\textbf{J}^{-1}_{SO(2)}(\mu)/SO(2)$.\\

In order to solve these two problems simultaneously, we will
define, in Section $2$, a precise notion of a non-smooth action of
a finite dimensional Lie group on a collection of smooth Banach
manifolds. We will see that the actions that appear in the
previous examples are particular cases of these non-smooth
actions. In Section $3$ we shall show that the non-smooth actions
we consider admit slices. This fact will be useful in the study
of the geometric properties of the orbit type sets, the isotropy type
sets, the fixed point sets, and of the orbit space. Due to the fact
that the action is not smooth, the isotropy strata are not smooth
manifolds. However we will prove in Section $4$ that, in the case
all the isotropy groups are conjugated, the orbit space is a
topological Banach manifold, by constructing explicit charts. Then
we shall prove that the changes of charts are $C^1$ with respect
to a weaker topology. This will allow us to define in Section $5$ a weak
tangent bundle for the orbit space, the notion of weak differentiable
curves, as well as the notion of differentiable functions on the
orbit space. The goal of Section $6$ is to apply
all the results of this paper to the case of the motion of the
incompressible fluid in order to carry out in a precise sense the
Poisson reduction by stages relative to the commuting actions $R$
and $L$ given in the first example.

\section{Non-smooth actions and their orbits}
\label{section: nonsmooth actions}

Let $\{Q^s|s> s_0\}$ be a collection of smooth Banach manifolds such that
for all $r>s> s_0$ there is a smooth inclusion $j_{(r,s)} :
Q^r\hookrightarrow Q^s$ with dense range satisfying  the following
condition: for all $q\in Q^r$, the range of the tangent map
$T_qj_{(r,s)}:T_q Q^r\longrightarrow T_q Q^s$ is dense. Density in these
two conditions  is always relative to the ambient spaces $Q^s$ and
$T_qQ^s$, respectively. In addition, we suppose that for each chart
$\varphi^s : U^s \longrightarrow T_qQ^s$ of $Q^s$ at $q\in Q^r,r>s$, the
map $\varphi^r:=\varphi^s|U^s\cap Q^r$ takes values in $T_qQ^r$ and is
a chart for $Q^r$. If these hypotheses hold, $\{Q^s\mid s >
s _0\} $ is called a \textbf{chain of Banach manifolds}.

Typical examples of chains of Banach manifolds are the collections of
manifolds of maps $Q^s:=H^s(M,N),
s>\frac{\operatorname{dim}M}{2}$, or
$Q^s:=C^s(M,N), s\geq 1$, where $M$ and $N$ are compact and
oriented finite dimensional manifolds ($M$ possibly with
boundary).

We suppose that each $Q^s,s>s_0$ carries a weak Riemannian metric
$\gamma_s$ which does not depend on $s$, that is
$j_{(r,s)}^*\gamma_s=\gamma_r$. So we will suppress the index $s$
and this metric will be denoted simply by $\gamma$.

For example, on $H^s(M,N)$ or $C^s(M,N)$ we can consider the $L^2$
metric
\[
\gamma(f)(u_f,v_f):=\int_Mg(f(x))(u_f(x),v_f(x))\mu(x)
\]
where $g$ is a Riemannian metric on $N$ and $\mu$ is a volume form
on $M$.

Let $G$ be a finite dimensional Lie group. We suppose that for all
$s>s_0$ we have a \textbf{proper} and \textbf{continuous} action
\[
\Phi : G\times Q^s \longrightarrow Q^s,\;\;\Phi(g,q)=\Phi_g
(q)=\Phi^q(g)
\]
of $G$ on $Q^s$ \textbf{compatible} with the restrictions $j_{(r,s)}^*$,
that is,
\begin{itemize}
\item $\Phi$ is a continuous map for every $s > s_0$;
\item $ \Phi$ is proper, which means that for each convergent sequences
$(q_n)_{n\in\mathbb{N}}$ and $(\Phi_{g_n}(q_n))_{n\in\mathbb{N}}$
in $Q^s$, there exists a convergent subsequence
$(g_{n_k})_{k\in\mathbb{N}}$ of $(g_n)_{n\in\mathbb{N}}$ in $G$;
\item each homeomorphism $\Phi_g$ commutes with all the maps $j_{(r,s)}$
which means that the action does not depend on $s $.

\end{itemize}
As a consequence of the properness assumption, we obtain that the
isotropy groups $G_q:=\{g\in G\,|\,\Phi_g(q)=q\}$ are compact.

Finally we suppose that:\\
$(1)$ For all $s> s_0$ and for all $g\in G$:
\begin{equation}
\label{phi_g_smooth} \Phi_g : Q^s\longrightarrow Q^s\text{ is
smooth.}
\end{equation}
$(2)$ For all $s> s_0+1$ and for all $q\in Q^s$:
\begin{equation}
\label{phi_c_one} \Phi^q : G\longrightarrow Q^{s-1}\text{ is }C^1.
\end{equation}
As a consequence of $(2)$, we obtain that for $s> s_0+1$ the
infinitesimal generators of this action are not vector fields on
$Q^s$ but they are sections of the smooth vector bundle
$TQ^{s-1}|Q^s\rightarrow Q^s$. Locally, this vector bundle is the
product of an open set in the model of $Q ^s$ and of the model
Banach space of $Q^{s-1} $. Indeed, denoting by $\xi_{Q^s}$ the
infinitesimal generator associated to $\xi\in\mathfrak{g}:=T_eG$,
we have, for all $q\in Q^s,s> s_0+1$:
\[
\xi_{Q^s}(q)=\left.\frac{d}{dt}\right|_{t=0}\Phi_{\operatorname{exp}(t\xi)}(q)=
T_e\Phi^q (\xi)\in T_qQ^{s-1}.
\]

\noindent \textbf{Example.} A typical example of such an action is
\[
\Phi:G\times H^s(G,M)\longrightarrow
H^s(G,M),\;\;\Phi(g,\eta):=\eta\circ R_g
\]
where $G$ is a compact finite dimensional Lie group, $M$ is a
compact finite dimensional manifold,
$s>\frac{\operatorname{dim}G}{2}$, and $R_g$ is the right
multiplication in $G$, that is, $R_g(h)=hg$. Using the fact that
$G$ can be viewed as a submanifold of $\mathcal{D}^r(G)$, the
group of class $H^r$-diffeomorphisms of $G$,
$r>\frac{\operatorname{dim} G}{2}+1$ (see Lemma \ref{submanifold}
below), and the fact that the composition
\[
\circ : \mathcal{D}^r(G)\times H^s(G,M)\longrightarrow
H^s(G,M),\;\;(\gamma,\eta)\longmapsto \eta\circ\gamma
\]
is continuous for $r\geq s$ (by Lemma $3.1$ of \cite{Eb1968}), we
obtain that $\Phi$ is a continuous action.

For all $g\in G$, the map
\[
\Phi_g : H^s(G,M)\longrightarrow
H^s(G,M),\;\;\Phi_g(\eta)=\eta\circ R_g
\]
is smooth by the $\alpha$-Lemma (see Proposition $3.4$ of
\cite{Eb1968}) and the tangent map is given by
$T\Phi_g(v_\eta)=v_\eta\circ R_g$, for $v_\eta\in T_\eta
H^s(G,M)$. So hypothesis $(\ref{phi_g_smooth})$ is verified.

We now check hypothesis $(\ref{phi_c_one})$. For all $\eta\in
H^s(G,M), s>\frac{\operatorname{dim}G}{2}+1$, the map
\[
\mathcal{D}^r(G)\longrightarrow H^{s-1}(G,M),\;\;\gamma\longmapsto
\eta\circ\gamma
\]
is of class $C^1$ by the proof of Proposition $3.4$ of
\cite{Eb1968}, for $r$ sufficiently large. Since $G$ is a
submanifold of $\mathcal{D}^r(G)$, the map
\[
\Phi^\eta : G\longrightarrow
H^{s-1}(G,M),\;\;\Phi^\eta(g)=\eta\circ R_g
\]
is of class $C^1$, and its tangent map is given by
$T\Phi^\eta(\xi_g)(h)=T\eta(TL_h(\xi_g))$, for $\xi_g\in T_gG$,
and where $L_g(h)=gh$. Remark that $T\Phi^\eta(\xi_g)$ is in
$T_{\Phi^\eta(g)}H^{s-1}(G,M)$ and does not belong to
$T_{\Phi^\eta(g)}H^s(G,M)$, in general. So the map
\[
\Phi^\eta : G\longrightarrow H^s(G,M)
\]
cannot be $C^1$ in general. \quad $\blacklozenge$

In Section 6 we will prove that with
$Q^s=\mathfrak{X}^s_{div}(M),
s>s_0=\frac{\operatorname{dim}M}{2}$, the action
\[
l:Iso^+\times \mathfrak{X}^s_{div}(M)\longrightarrow
\mathfrak{X}^s_{div}(M),
\]
defined in Section $1$, is continuous and verifies the hypotheses
(\ref{phi_g_smooth}) and (\ref{phi_c_one}).

In a similar way, for the second example treated in the Section 1,
choosing $Q^s=P^s/SO(2)$, $s>s_0=\frac{3}{2}$, we see that the
action of $S^1$ induced by $R$ is continuous and verifies the same
hypotheses.

\medskip

We now prove the following Lemma.

\begin{lemma} \label{submanifold} $G$ is a submanifold of $\mathcal{D}^r(G)$,
for all $r>\frac{\operatorname{dim} G}{2}+1$.
\end{lemma}
\textbf{Proof.} We identify $G$ with the subgroup $\{R_g :
G\longrightarrow G\,|\,g\in G\}$ of $\mathcal{D}^r(G)$. Recall
that the model of $\mathcal{D}^r(G)$ is the Hilbert space
$T_e\mathcal{D}^r(G)=\mathfrak{X}^r(G)$, consisting of the class
$H^r$ vector fields on $G$. A chart of $\mathcal{D}^r(G)$ at
$e=id_G$ is given by
\[
\psi : U\subset
\mathfrak{X}^r(G)\longrightarrow\mathcal{D}^r(G),\;\;\psi(X)=\operatorname{Exp}
\circ X,
\]
where $\operatorname{Exp} : TG\longrightarrow G$ is the map defined
by  $\operatorname{Exp}|T_gG=\operatorname{exp}_g$, for
$\operatorname{exp}_g(\xi_g):=L_g(\operatorname{exp}(TL_{g^{-1}}(\xi_g)))$,
and $\exp: \mathfrak{g} \rightarrow G $ the exponential map of the Lie
group $G$.

Let
$\mathfrak{X}_L(G):=\{X\in\mathfrak{X}(G)\,|\,(L_g)^*X=X,\forall\,g\in
G\}$ be the space of left-invariant vector fields on $G$. Then we
have the isomorphism $\mathfrak{g}\longrightarrow
\mathfrak{X}_L(G),\;\;\xi\longmapsto X_\xi^L$ where
$X_\xi^L(g):=T_eL_g(\xi)$. Since $\mathfrak{X}_L(G)$  is finite
dimensional, it is a closed subspace of the Hilbert space
$\mathfrak{X}^r(G)$ and the topology induced by
$\mathfrak{X}^r(G)$ on $\mathfrak{X}_L(G)$ coincides with the one
induced by the identification with $\mathfrak{g}$. To show
that $\mathfrak{X}_L(G)$ is in fact the model of $G$ viewed as a
submanifold of $\mathcal{D}^r(G)$, it suffices to see that
$\psi(U\cap\mathfrak{X}_L(G))=\psi(U)\cap G.\;\;\;\; \blacksquare$

\medskip

Denote by $\operatorname{Orb}(q):=\{\Phi_g(q)\,|\,g\in G\}\subset
Q^s$ the orbit of $q\in Q^s$. Since $\xi_{Q^s}(q)$ is not in
$T_qQ^s$, $\operatorname{Orb}(q)$ can not be a submanifold of
$Q^s$. However we have the following result.

\begin{theorem} For all $q\in Q^s, s>s_0+1$, $\operatorname{Orb}(q)$
is a submanifold of $Q^{s-1}$ and
$T_q\operatorname{Orb}(q)=\{\xi_{Q^s}(q)\,|\,\xi\in\mathfrak{g}\}\subset
T_qQ^{s-1}$.
\end{theorem}
\textbf{Proof.} The $C^1$ map $\Phi^q : G\longrightarrow
Q^{s-1}$ induces the $C^1$ map $\widetilde{\Phi}^q :
G/G_q\longrightarrow Q^{s-1}$ defined by
$\Phi^q=\widetilde{\Phi}^q\circ \pi_{G,G_q}$, where $\pi_{G,G_q} :
G\longrightarrow G/G_q$ is the projection given by
$\pi_{G,G_q}(g)=gG_q$. Since
$\widetilde{\Phi}^q(G/G_q)=\operatorname{Orb}(q)$, we will prove
that $\widetilde{\Phi}^q$ is an embedding. It suffices to show
that $\widetilde{\Phi}^q$ is a closed injective immersion. This
can be proven like in the usual case of a smooth action on a
finite dimensional manifold (see Corollary 4.1.22 in
\cite{AbMa1978} for example). It just remains to show that the
range of $T_{[e]}\tilde{\Phi}^q$ is closed and split in
$T_qQ^{s-1}$. This is true because
$T_{[e]}\widetilde{\Phi}^q(T_{[e]}(G/G_q))$ is a finite
dimensional vector space.$\;\;\;\; \blacksquare$

\section{The slice theorem and its consequences}

Recall that in the case of a smooth and proper action $\Phi :
G\times M\longrightarrow M$ of a finite dimensional Lie group $G$
on a finite dimensional manifold $M$, the existence of slices is
a key fact in the the study of the geometric properties of the
submanifolds $M^H$, $M_H$, and $M_{(H)}$ and of the
stratification of the orbit space $M/G$. In the more general
case of a continuous action on an infinite dimensional Banach
manifold, we adopt the following definition of a slice (see for
example Theorem 7.1 in
\cite{Eb1968} or Theorem 4.1 in \cite{IsMa1982}).

\begin{definition} A \textbf{slice} at $q\in Q^s$ is a submanifold
$S_q\subset Q^s$ containing $q$ such that:\\
$(\rm{S}1)$ if $g\in G_q$, then $\Phi_g(S_q)=S_q$,\\
$(\rm{S}2)$ if $g\in G$ and $\Phi_g(S_q)\cap S_q\neq\varnothing$,
then
$g\in G_q$,\\
$(\rm{S}3)$ there is a local section $\chi : G/G_q\longrightarrow
G$ defined in a neighborhood $V([e])$ of $[e]$ in $G/G_q$ such
that the map
\[
F : V([e])\times S_q\longrightarrow
Q^s,\;\;F([g],q):=\Phi_{\chi([g])^{-1}}(q)
\]
is an homeomorphism onto a neighborhood $U$ of $q$.
\end{definition}

Note that usually we take $F([g],q):=\Phi_{\chi([g])}(q)$ but it
will be more natural to use $F([g],q):=\Phi_{\chi([g])^{-1}}(q)$.

In finite dimensions and in the case of a smooth action, the
previous definition is equivalent to the usual ones (see
\cite{OrRa2004} for equivalent definitions of a slice). We will
show that slices exist for a non-smooth action $\Phi$ verifying
all the hypotheses in Section 2.

Note that in \cite{Pa1961} it is proven that slices exist for each
continuous and proper (in the sense given there) action $G\times
X\longrightarrow X$ of a finite dimensional Lie group $G$ on a
completely regular topological space $X$; nevertheless, the
construction is not explicit and so it can not be used below to
construct charts for $Q^s/G$.

Recall that the action we consider is proper, so the isotropy
groups $G_q$ are compact. Now we prove that for each $q\in Q^s$,
there exist a $G_q$-invariant chart of $Q^s$ at $q$. This result
is shown, for example, in the Appendix B of \cite{CuBa1997}, in
the case of a smooth and proper action on a finite dimensional
manifold. We give below the proof of this statement in order to
show carefully that it is still valid in our case of a non-smooth
action on an infinite dimensional Banach manifold.

\begin{lemma}\label{H-invariant charts}
Assume that the action $\Phi : G\times
Q^s\longrightarrow Q^s$
verifies all the hypotheses  in  Section 2. Then for all
$q\in Q^s$ there
exists a chart $(\psi,V)$ of $Q^s$ at $q$ such that:\\
$(1)$ $V$ is $G_q$-invariant, and $\psi : V\subset
Q^s\longrightarrow
T_qQ^s$ verifies $\psi(q)=0$ and $T_q\psi=\operatorname{id}_{T_qQ^s}$,\\
$(2)$ $\forall\,r\in V,\forall\,h\in G_q$, we have:
\[
\psi(\Phi_h(r))=T_q\Phi_h(\psi(r)).
\]
\end{lemma}
\textbf{Proof.} Let $(\varphi,U)$ be a chart of $Q^s$ at $q$ such
that $\varphi(q)=0_q$ and $T_q\varphi=\operatorname{id}_{T_qQ^s}$.
Since $G_q$ is compact, there exists a $G_q$-invariant
neighborhood $V\subset U$ of $q$; see for example Lemma 2.3.29 of
\cite{OrRa2004} which does not use in its proof the finite
dimensionality of the manifold. Let
$W:=\varphi(V)$ and let
$\overline{\Phi}$ be the action of $G_q$ induced on $W$ by
$\varphi$, that is, we have the commuting diagram
$$\xymatrix{
\Phi : G_q\times V\subset Q^s\;\;\;\ \ar[r] \ar[d]_{\varphi}&
V\ar [d]^{\varphi}\\
\overline{\Phi} : G_q\times W\subset T_qQ^s \ar[r] & W. }$$ We
consider the map $\overline{\psi} : W\subset T_qQ^s\longrightarrow
T_qQ^s$ defined by
\[
\overline{\psi}(v_q):=\int_{G_q}D\overline{\Phi}_{g^{-1}}(0)(\overline{\Phi}_g
(v_q))dg
\]
where $D\overline{\Phi}_{g^{-1}}$ is the Fr\'echet derivative of
the smooth map $\overline{\Phi}_{g^{-1}}:W\subset
T_qQ^s\longrightarrow T_qQ^s$ and $dg$ is the Haar measure of the
compact group $G_q$ such that $\operatorname{Vol}(G_q)=1$.

With this definition of $\overline{\psi}$ we obtain that (see p.
301 of \cite{CuBa1997} for details):
\[
D\overline{\Phi}_h(0)(\overline{\psi}(v_q))=\overline{\psi}(\overline{\Phi}_h
(v_q)),
\forall\, h\in G_q, \forall\, v_q\in W\;\text{ and
}\;D\overline{\psi}(0)=\operatorname{id}_{T_qQ^s}.
\]
Define $\psi:=\overline{\psi}\circ\varphi : V \longrightarrow
T_qQ^s$. Then  $T_q\psi=\operatorname{id}_{T_qQ^s}$ and we have
\begin{align*}
&D\overline{\Phi}_h(0)(\overline{\psi}(v_q))=\overline{\psi}(\overline{\Phi}
(v_q))\\
&\Longrightarrow
D(\varphi\circ\Phi_h\circ\varphi^{-1})(0)(\psi(\varphi^{-1}(v_q)))=\psi(\Phi_h
(\varphi^{-1}(v_q)))\\
&\Longrightarrow
T_q\varphi(T_q\Phi_h(T_0\varphi^{-1}(\psi(\varphi^{-1}(v_q)))))=\psi(\Phi_h
(\varphi^{-1}(v_q)))\\
&\Longrightarrow
T_q\Phi_h(\psi(\varphi^{-1}(v_q)))=\psi(\Phi_h(\varphi^{-1}(v_q))), \text{
since $T_q\varphi=\operatorname{id}_{T_qQ^s}$}\\
&\Longrightarrow T_q\Phi_h(\psi(r))=\psi(\Phi_h(r)), \forall \,
r\in V. \;\;\;\; \blacksquare
\end{align*}

We will define below a map $\mathcal{S}_q$  that will play a
central role in the rest of the paper and is a generalization of
the map
$S$ defined in Chapter 7 of \cite{Karcher1970} in the case of the
natural $S^1$-action on the Hilbert manifold $H^1(S^1,M)$ of
closed $H^1$-curves in a compact Riemannian manifold $M$ without
boundary. In order to define this map $\mathcal{S}_q$ we will need
the following Lemma.

\begin{lemma} Let $\Phi :G\times Q^s\longrightarrow Q^s$ be a
continuous action. Let $q\in Q^s$ be such that $H:=G_q$ is
compact and let $U_1$ be a neighborhood of $q$ in $Q^s$. Then
there exist a neighborhood $V(H)$ of $H$ in $G$ and a
neighborhood $U_2$ of
$q$ in $Q^s$ such that:
\[
\Phi(V(H)\times U_2)\subset U_1.
\]
\end{lemma}
\textbf{Proof.} We have $\Phi(h,q)=q, \forall\,h\in H$, so for all
$h\in H$, there exist a neighborhood $V_h$ of $h$ in $H$ and a
neighborhood $U_2^h$ of $q$ in $Q^s$ such that
\[
\Phi(V_h\times U_2^h)\subset U_1.
\]
Since $(V_h)_{h\in H}$ is an open covering of the compact
group $H$, there exists a finite subcovering
$(V_{h_i})_{i=1...n}$. Let
\[
V(H):=\bigcup_{i=1}^n V_{h_i} \text{ and } U_2:=\bigcap_{i=1}^n
U^{h_i}_2.
\]
Clearly $V(H)$ and $U_2$ are neighborhoods of $H$ and $q$;
moreover:
\[
\Phi(V(H)\times U_2)=\Phi\Big{(}\bigcup_{i=1}^n V_{h_i}\times
U_2\Big{)}\subset\bigcup_{i=1}^n \Phi(V_{h_i}\times
U_2)\subset\bigcup_{i=1}^n \Phi(V_{h_i}\times U_2^{h_i})\subset
U_1.\;\;\;\; \blacksquare
\]

Recall that $G/H$ is a smooth manifold (not a group in general)
whose elements will be denoted by $gH$ or $[g]$;
$\pi_{G,H} : G\longrightarrow G/H$ is the projection map.

\medskip

Assume now that the action $\Phi : G\times
Q^s\longrightarrow Q^s$
verifies all the hypotheses  in  Section 2. Let $q\in Q^s,
s>s_0+1$, $H:=G_q$, and $\mathfrak{h}:=T_eH$. Let
$\operatorname{B}:=(e_1,...,e_k,e_{k+1},...,e_n)$ be a basis of
$\mathfrak{g}$ such that $(e_1,...,e_k)$ is a basis of
$\mathfrak{h}$. Let $(\varphi^{s-1},U^{s-1})$ be a chart of
$Q^{s-1}$ at $q\in Q^s$ such that $\varphi^{s-1}(q)=0_q$ and
$T_q\varphi^{s-1}=\operatorname{id}_{T_qQ^{s-1}}$. Let $\chi :
G/H\longrightarrow G$ be a local section defined in a neighborhood
of $[e]$. We define the map
\[
\mathcal{S}_q : V([e])\times U\subset G/H\times Q^s\longrightarrow
\mathfrak{n},\;\;\mathcal{S}_q(gH,r):=\sum_{i=k+1}^n\gamma(q)\left(\varphi^{s-
1}(\Phi_{\chi(gH)}(r)),E_i(q)\right)e_i
\]
where the neighborhoods $V(H)$ and $U$ are such that
$\Phi(V(H)\times U)\subset U^{s-1}\cap Q^s$ (which is possible by
the preceding Lemma), $V([e])$ is such that $\chi(V([e]))\subset
V(H)$, $\mathfrak{n}$ is the subspace of $\mathfrak{g}$ generated
by $(e_{k+1},...,e_n)$, and $E_i(q):=T_e\Phi^q(e_i)\in T_qQ^{s-1}$
is the infinitesimal generator associated to $e_i$.

Remark that the map $\mathcal{S}_q$ depends on the basis
$\operatorname{B}$ and on the local section $\chi$. Note that
$\mathcal{S}_q(eH,q)=0$. Note also that if $\Phi$ is
free at $q$, we have $H=\{e\}$ and the map $\mathcal{S}_q$ is
given by
\begin{equation}
\label{s_q_useful_formula}
\mathcal{S}_q : V(e)\times U\subset G\times Q^s \longrightarrow
\mathfrak{g},\;\;\mathcal{S}_q(g,r):=\sum_{i=1}^n\gamma(q)\left(\varphi^{s-1}
(\Phi_g(r)),E_i(q)\right)e_i.
\end{equation}
The following Lemma states the most important property of the map
$\mathcal{S}_q$.

\begin{lemma} The map $\mathcal{S}_q$ is of class $C^1$.
Moreover, if we consider the map
$\mathcal{S}_q(\_,q) :V([e])\subset G/H\longrightarrow
\mathfrak{n}$, then its tangent map at $eH$,
\[
T_{eH}[\mathcal{S}_q(\_,q)]: T_{eH}(G/H)\longrightarrow
\mathfrak{n}
\]
is invertible. More precisely, its matrix relative to the basis
$(e_{k+1},...,e_n)$ is given by:
\[
\Big{(}T_{eH}[\mathcal{S}_q(\_,q)]\Big{)}=\Big{(}\gamma(q)(E_i(q),E_j(q))\Big
{)}_{i,j=k+1,...,n}
\]

\end{lemma}
\textbf{Proof.} By the assumptions (\ref{phi_g_smooth}) and
(\ref{phi_c_one}) on the action, we obtain that the map
\[V([e])\times U \longrightarrow U^{s-1},\;\;
(gH,r)\longmapsto \Phi_{\chi(gH)}(r)
\]
is $C^1$ as a map with values in $Q^{s-1}$. Using that
$\varphi^{s-1}$ is a smooth chart for $Q^{s-1}$ and that the
bilinear form $\gamma(q)$ is continuous, we obtain that
$\mathcal{S}_q$ is $C^1$. Now we compute the tangent map to
$\mathcal{S}_q(\_,q) :V([e])\subset G/H\longrightarrow
\mathfrak{n}$ at $eH$. Consider the tangent map $T_e\pi_{G,H} :
\mathfrak{g}\longrightarrow T_{eH}(G/H)$; clearly we have
$\operatorname{ker}(T_e\pi_{G,H})=\mathfrak{h}$ and we obtain that
$T_e\pi_{G,H} : \mathfrak{n}\longrightarrow T_{eH}(G/H)$ is
bijective, so $(T_e\pi_{G,H}(e_{k+1}),...,T_e\pi_{G,H}(e_n))$ is a
basis of $T_{eH}(G/H)$. Since
$\pi_{G,H}\circ\operatorname{exp}(te_j)$ is a curve in $G/H$
tangent to $T\pi_{G,H}(e_j)$ at $eH$, we have:
\begin{align*}
T_{eH}[\mathcal{S}_q(\_,q)](T\pi_{G,H}(e_j))&=\left.\frac{d}{dt}\right|_{t=0}
\mathcal{S}_q(\pi_{G,H}(\operatorname{exp}(te_j)),q)\\
&=\left.\frac{d}{dt}\right|_{t=0}\sum_{i=k+1}^n\gamma(q)\left(\varphi^{s-1}
(\Phi_{\chi(\pi_{G,H}(\operatorname{exp}(te_j)))}(q)),E_i(q)\right)e_i\\
&=\left.\frac{d}{dt}\right|_{t=0}\sum_{i=k+1}^n\gamma(q)\left(\varphi^{s-1}
(\Phi_{\operatorname{exp}(te_j)}(q)),E_i(q)\right)e_i\\
&=\sum_{i=k+1}^n\gamma(q)\left(T_q\varphi^{s-1}\left(\left.\frac{d}{dt}\right|_
{t=0}\Phi_{\operatorname{exp}(te_j)}(q)\right),E_i(q)\right)e_i\\
&=\sum_{i=k+1}^n\gamma(q)(E_j(q),E_i(q))e_i.
\end{align*}
So we obtain that
\[
\Big{(}T_{eH}[\mathcal{S}_q(\_,q)]\Big{)}=\Big{(}\gamma(q)(E_i(q),E_j(q))\Big
{)}_{i,j=k+1,...,n}
\]
and since $\gamma(q)$ is strongly non-degenerate on any finite
dimensional vector space, it remains to show that
$(E_{k+1}(q),...,E_n(q))$ is linearly independent. To prove this,
we consider the tangent map $T_e\Phi^q :
\mathfrak{g}\longrightarrow T_qQ^{s-1}$. We have
$\operatorname{ker}(T_e\Phi^q)=\mathfrak{h}$. Indeed, if
$\xi\in\mathfrak{g}$ is such that $T_e\Phi^q(\xi)=0$ then we
have:
\begin{align*}
\frac{d}{dt}\Phi_{\operatorname{exp(t\xi)}}(q)&=\left.\frac{d}{ds}\right|_{s=0}
\Phi_{\operatorname{exp((t+s)\xi)}}(q)=\left.\frac{d}{ds}\right|_{s=0}\Phi_
{\operatorname{exp(t\xi)}}(\Phi_{\operatorname{exp(s\xi)}}(q))\\
&=T_q\Phi_{\operatorname{exp(t\xi)}}(T_e\Phi^q(\xi))=0,
\end{align*}
and we conclude that
$\Phi_{\operatorname{exp(t\xi)}}(q)=\Phi_{\operatorname{exp(0\xi)}}(q)=q$
for all $t$, so $\operatorname{exp(t\xi)}\in G_q=H$ and
$\xi\in\mathfrak{h}$. Thus we obtain that $T_e\Phi^q :
\mathfrak{n}\longrightarrow T_qQ^{s-1}$ is injective, and
$(E_{k+1}(q),...,E_n(q))$ is linearly independent since it is the
image of the basis $(e_{k+1},...,e_n)$ of $\mathfrak{n}$ under the
injective linear map $T_e\Phi^q$. $\;\;\;\; \blacksquare$

\medskip

For $q\in Q^s, s>s_0+1$, we define the closed codimension $(n-k)$
subspace
\begin{equation}
\label{n q s def}
N^s_q:=\{v_q\in T_qQ^s\,|\,\gamma(q)(v_q,E_i(q))=0, \text{ for all
} i=k+1,...,n\}
\end{equation}
in $T_qQ^s$, where as before,
$\operatorname{B}=(e_1,...,e_k,e_{k+1},...,e_n)$ is a basis of
$\mathfrak{g}$ such that $(e_1,...,e_k)$ is a basis of
$\mathfrak{h}:=T_eG_q$, and $E_i(q)$ are the corresponding
infinitesimal generators. We can now state the following useful
consequence of the properties of $\mathcal{S}_q$.

\begin{theorem}\label{IFThm} There exists a neighborhood $V([e])$
of $eH$ in $G/H$, a neighborhood $U$ of $q$ in $Q^s$, and a map
$\beta_q : U\longrightarrow V([e])$ of class $C^1$ such that for
all $[g]\in V([e])$ and $r\in U$ we have:
\[
\varphi^{s-1}(\Phi_{\chi([g])}(r))\in N^s_q \Longleftrightarrow
\mathcal{S}_q([g],r)=0 \Longleftrightarrow [g]=\beta_q(r).
\]
\end{theorem}
\textbf{Proof.} The first equivalence follows from the definition
of $\mathcal{S}_q$. For the second it suffices to use the implicit
function theorem since the map $\mathcal{S}_q$ is of class $C^1$
and verifies $\mathcal{S}_q(eH,q)=0$ and
$T_{eH}[\mathcal{S}_q(\_,q)]$ is invertible, by the preceding
Lemma.$\;\;\;\; \blacksquare$

\begin{theorem} (Slice theorem) Assume that the action verifies
all the hypotheses in Section 2. Then for each $q\in Q^s,
s>s_0+1$, there exists a slice at $q$.
\end{theorem}
\textbf{Proof.} Let $H:=G_q$. By Lemma \ref{H-invariant charts},
there exists a $H$-invariant chart $(\psi,V)$ of $q$ at $Q^s$.
Let $\overline{\gamma}$ be the metric defined by
\[
\overline{\gamma}(q)(u_q,v_q):=\int_H\gamma(\Phi_h(q))(T\Phi_h(u_q),T\Phi_h
(v_q))dh,
\]
where $dh$ is the Haar measure on $H$ such that
$\operatorname{Vol}(H)=1$. Then we have
\[
(\Phi_h)^*\overline{\gamma}=\overline{\gamma},\forall\,h\in H.
\]
Consider the closed codimension $(n-k)$ subspace $N^s_q:=\{v_q\in
T_qQ^s\,|\,\overline{\gamma}(q)(v_q,E_i(q))=0, \text{ for all }
i=k+1,...,n\}$, and define the submanifold
\[
\widetilde{S}_q:=\psi^{-1}(N^s_q\cap \psi(V)).
\]
By invariance of $\psi$ and $N^s_q$, we obtain that
$\Phi_h(\widetilde{S}_q)=\widetilde{S}_q$ for all $h\in H$, so
condition $(\rm{S}1)$ is satisfied.

We now check condition $(\rm{S}3)$. Let $\chi : W([e])\subset
G/H\longrightarrow G$ be a local section defined in a neighborhood
$W([e])$ of $[e]$ and consider the continuous map
\[
F : V([e])\times S_q\longrightarrow
Q^s,\;\;F([g],q):=\Phi_{\chi([g])^{-1}}(q).
\]
By Theorem \ref{IFThm}, there is a neighborhood
$V([e])$ of $eH$ in $G/H$, a neighborhood $U$ of $q$ in $Q^s$, and
a $C^1$ map $\beta_q : U\longrightarrow V([e])$ such that
\[
\psi(\Phi_{\chi([g])}(r))\in N^s_q \Longleftrightarrow
[g]=\beta_q(r).
\]
By definition of $\widetilde{S}_q$ we obtain
\[
\Phi_{\chi([g])}(r)\in\widetilde{S}_q\Longleftrightarrow
[g]=\beta_q(r),
\]
so we can define the continuous map
\[
U\subset Q^s\longrightarrow
V([e])\times\widetilde{S}_q,\;\;r\longmapsto
(\beta_q(r),\Phi_{\chi(\beta_q(r))}(r))
\]
and one sees that it is the inverse of the map $F$. This proves
that $F$ is an homeomorphism. Remark that for $(\rm{S}1)$
and$(\rm{S}3)$ we do not use the properness of the action.

We now show that there exists a neighborhood $U$ of $q$ such that
choosing $S_q:=\widetilde{S}_q\cap U$ we have condition
$(\rm{S}2)$:
\[
\Phi_g(S_q)\cap S_q\neq \varnothing\Longrightarrow g\in H.
\]
By contradiction, suppose that such a neighborhood does not exist.
So, without loss of generality, we can consider a $H$-invariant
fundamental system of neighborhoods $(U^m)_{m\in\mathbb{N}}$ of
$q$ such that for all $m\in\mathbb{N}$ we have the following property:
\begin{quote}{there exists $g_m\notin H$ such that
$\Phi_{g_m}(S^m_q)\cap S^m_q\neq\varnothing$, where
$S^m_q:=\widetilde{S}_q\cap U^m$.}\end{quote}

Note that $\Phi_{g_m}(S^m_q)\cap S^m_q\neq\varnothing$ if and
only if $\Phi_{g_m ^{-1}}(S^m_q)\cap S^m_q\neq\varnothing$.
So for all $m$ there exists $g_m\notin H$ and $r_m\in S^m_q$ such
that $\Phi_{g_m}(r_m)\in S^m_q$. Then we have
$r_m\longrightarrow q$ and $\Phi_{g_m}(r_m)\longrightarrow q$
and by properness of the action we obtain the existence of a subsequence
$(g_{m_k})_{k\in\mathbb{N}}\subset (g_m)_{m\in\mathbb{N}}$ and of
an element $g\in G$ such that $g_{m_k}\longrightarrow g$. Therefore
$\Phi_g(q)=q$ and hence $g\in H$.

Let $V([e])$ be a neighborhood of $[e]$ in $G/H$ and
$V(H):=\pi_{G,H}^{-1}(V([e])$ a neighborhood of $H$ in $G$. We can
suppose that $g_{m_k}\in V(H)$. Since $\Phi_{g_{m_k}}(r_{m_k})\in
S^{m_k}_q$, we have
\[
S^{m_k}_q \ni
\Phi_{g_{m_k}}(r_{m_k})=\Phi_{\chi(\pi_{G,H}(g_{m_k}))h}(r_{m_k})=\Phi_{\chi
(\pi_{G,H}(g_{m_k}))}(\Phi_h(r_{m_k})),
\]
where $h\in H$. For sufficiently small $V([e])$, we have
$[g_{m_k}]=\beta_q(\Phi_h(r_{m_k}))$. But by $(\rm{S}1)$ we know
that $\Phi_h(r_{m_k})\in S^m_q$ and hence
$\beta_q(\Phi_h(r_{m_k}))=eH$. Thus we obtain that $[g_{m_k}]=eH$
and so $g_{m_k}\in H$. This last affirmation is a
contradiction.$\;\;\;\; \blacksquare$

\medskip

Here are two immediate consequences of the Slice Theorem:\\
$(1)$ For all $s\in F(V([e])\times S_q)\cap S_q$ we have
$G_s\subset G_q$.\\
Indeed, $g\in G_s$ implies $\Phi_g(s)=s$, but since $s\in S_q$ we
obtain $\Phi_g(S_q)\cap S_q\neq \varnothing$. So by $(\rm{S}2)$
we have $g\in G_q$.\\
$(2)$ For all $r\in F(V([e])\times S_q)$, $G_r$ is conjugated to a
subgroup of $G_q$.\\
Indeed, $r\in F(V([e])\times S_q)$ implies
$r=\Phi_{\chi([g])^{-1}}(s)$ with $s\in S_q$. So $G_r$ is
conjugated to $G_s$. Since $s\in F(V([e])\times S_q)\cap S_q$, by
$(1)$ we obtain that $G_s\subset G_q$.

\medskip

We now study the geometric properties of some important subsets of
$Q^s$. If $H$ is a subgroup of $G$, we denote by
\[
(Q^s)^H:=\{q\in Q^s\,|\,G_q \supset H\}
\]
the \textbf{$H$-fixed point set}, by
\[
(Q^s)_H=\{q\in Q^s\,|\,G_q=H\}
\]
the \textbf{$H$-isotropy type set}, and by
\[
(Q^s)_{(H)}:=\{q\in Q^s\,|\, G_q \text{ is conjugated to } H\}
\]
the \textbf{$(H)$-orbit type set}. We have $(Q^s)_H\subset
(Q^s)^H$. It is well known that in the case of a smooth and proper
action $\Phi : G\times M\longrightarrow M$ on a finite dimensional
manifold $M$, we obtain that $M^H, M_H$, and $M_{(H)}$ are
submanifolds of $M$ and that $M_H$ is an open subset of $M^H$. We
will see that some of these results are still valid in our case.

Define for $q \in (Q ^s)^H $ the closed subspace
\[
[T_qQ^s]^H:=\{v_q\in T_qQ^s\,|\,T_q\Phi_h(v_q)=v_q, \forall\,h\in
H \}.
\]

\begin{theorem}  Assume that the action $\Phi : G\times Q^s\longrightarrow Q^s$
verifies all the hypotheses in Section 2. Let $H$ be a
compact subgroup of $G$. Then $(Q^s)^H$ is a  smooth
submanifold of $Q^s$ whose tangent space at $q\in (Q^s)^H$ is
$T_q(Q^s)^H=[T_qQ^s]^H$.
\end{theorem}
\textbf{Proof.} If $q\in (Q^s)^H$, using the preceding lemma and
the fact that $H\subset G_q$, there exists a chart $(\psi,V)$ of
$Q^s$ at $q$ such that $\psi(\Phi_h(r))=T_q\Phi_h(\psi(r)),
\forall\,r\in V,\forall\,h\in G_q$. So we obtain that
$\psi(V\cap (Q^s)^H)=\psi(V)\cap [T_qQ^s]^H$, which proves that
$(Q^s)^H$ is a submanifold of $Q^s$ whose tangent space at $q $ is given
by $T_q(Q^s)^H=[T_qQ^s]^H.\;\;\;\; \blacksquare$

\begin{theorem}
\label{open} Assume that the action $\Phi : G\times
Q^s\longrightarrow Q^s$ verifies all the hypotheses  in
Section 2. Let $H$ be a compact subgroup of $G$. Then for all
$s>s_0+1$, $(Q^s)_H$ is an open subset of the manifold $(Q^s)^H$.
\end{theorem}
\textbf{Proof.} Let $q\in (Q^s)_H$. We will find a neighborhood of
$q$ in $(Q^s)^H$ which is included in $(Q^s)_H$. It suffices to
consider the open set $(Q^s)^H\cap F(V([e])\times S_q)$ in
$(Q^s)^H$. Indeed for each $r\in (Q^s)^H\cap F(V([e])\times S_q)$
we have $H\subset G_r$ and by the consequence $(2)$ of the Slice
Theorem, we know that $G_r$ is conjugated to a subgroup of
$G_q=H$. So we must have $G_r=H$ and we obtain $r\in (Q^s)_H$.
This proves that $(Q^s)^H\cap F(V([e])\times S_q)$ is an open
neighborhood of $q$ in $(Q^s)_H$.$\;\;\;\; \blacksquare$

\begin{corollary} Assume that the action $\Phi : G\times
Q^s\longrightarrow Q^s$ verifies all the hypotheses  in
Section 2. Let $H$ be a compact subgroup of $G$. Then for
$s>s_0+1$, $(Q^s)_H$ is a submanifold of $Q^s$.
\end{corollary}

We have not succeeded in proving that the $(H)$-orbit type set
$(Q^s)_{(H)}$ is a submanifold of $Q^s$. In fact, all the proofs
of this statement in the usual finite dimensional situation use the
smoothness of $\Phi$ on $G\times M$ in a crucial manner. However,
we will show in the next section in which sense $(Q^s)_{(H)}$
carries a smooth structure.

\section{Geometric properties of the orbit space and of the $(H)$-orbit type
sets}

In this section we assume that the action $\Phi$ verifies all the
hypotheses made in Section 2. Moreover we suppose that
\textbf{all the isotropy groups are conjugated to a given one},
say $H$. We consider the \textbf{orbit space} $Q^s/G$ endowed with
the quotient topology. Recall that the projection $\pi :
Q^s\longrightarrow Q^s/G$ is a continuous open map. We will denote
by $\omega,\upsilon...\in Q^s/G$ the orbits and by
$q_\omega,q_\upsilon,...\in Q^s$ the elements of $Q^s$ that belong
to the corresponding orbits, that is,  $q_\omega\in
\pi^{-1}(\omega),q_\upsilon\in \pi^{-1}(\upsilon)$.

The first goal of this section is to show that the orbit space
$Q^s/G$ is a (Banach) topological manifold. In order to do that,
we will construct explicit charts for $Q^s/G$. We will see that
$Q^s/G$ is not a smooth manifold. However, with these charts, we
will show in the second part of this section, that $Q^s/G$
shares many properties with a $C^1$-manifold. For example, we
will be able to define a notion of tangent bundle as well as a notion
of differentiable function for $Q^s/G$. Our construction of charts
is inspired by the method used by H. Karcher to prove the local
contractibility of the set $H^1(S^1,M)/S^1$ (see Chapter 7 of
\cite{Karcher1970}). This construction is done in several steps
given in the  following two lemmas. Recall at this point the defintion
of $N_q^s$ in \eqref{n q s def}.

\begin{lemma}
\label{A} Let $q\in Q^s, s> s_0+1$, and let $(\psi,V)$ be a
$G_q$-invariant chart of $Q^s$ at $q$ (whose existence is proved in
Lemma \ref{H-invariant charts}). Then for
$V$ sufficiently small, the continuous map
\[
\pi\circ \psi^{-1} : N_q^s\cap\psi(V) \longrightarrow Q^s/G
\]
is injective.
\end{lemma}
\textbf{Proof.} For $V$ sufficiently small,
$S_q:=\psi^{-1}(N^s_q\cap\psi(V))$ is a slice at $q$. Let
$u_q,v_q\in N^s_q\cap\psi(V)$ such that
$\pi(\psi^{-1}(u_q))=\pi(\psi^{-1}(v_q))$. Then we have
$\pi(s_1)=\pi(s_2)$ with $s_1=\psi^{-1}(u_q),
s_2=\psi^{-1}(v_q)\in S_q$. So we have
$s_1=\Phi_g(s_2)$ and we obtain $\Phi_q(S_q)\cap
S_q\neq\varnothing$. By $(\rm{S}2)$, we have $g\in G_q$. By the
consequence $(1)$ of the Slice Theorem, we know that
$G_{s_2}\subset G_q$, but since all the isotropy groups are
conjugated, we have $G_{s_2}=G_q$. So we obtain that
$s_1=\Phi_g(s_2)=s_2$. This proves that
$u_q=v_q$.$\;\;\;\;\blacksquare$

\begin{lemma}
\label{B} Let $q\in Q^s, s> s_0+1$, and let $(\psi,V)$ be a chart
of $Q^s$ at $q$. We can choose a neighborhood $U$ of $q$ in $Q^s$
such that the continuous map
\[
B_q : U \subset Q^s \longrightarrow N^s_q,\;\;B_q
(r):=\psi(\Phi_{\chi(\beta_q(r))}(r))
\]
does not depend on the choice of the representative in the equivalence
class of $r$ intersected with $U$. So, for all $\omega\in Q^s/G$
and $q_\omega\in\pi^{-1}(\omega)$ we can define the map $\mathcal
{B}_{q_\omega}$ such that the following diagram commutes:
\[
\xymatrix{
    U \ar[d]_{\pi} \ar[rd]^{B_{q_\omega}} & \\
    \mathcal{U} \ar[r]_{\mathcal{B}_{q_\omega}} &
N^s_{q_\omega}, }
\]
where $\mathcal{U}:=\pi(U)$ is a neighborhood of $\omega$ in
$Q^s/G$. Furthermore, the map $\mathcal{B}_{q_\omega}$ is
continuous and injective.
\end{lemma}
\textbf{Proof.} Let $r,\bar{r}\in U$ be in the same orbit, that
is, $\bar{r} =\Phi_g(r)$. We show that
$B_{q_\omega}(r)=B_{q_\omega}(\bar{r})$. We have
\[
\Phi_{\chi(\beta_q(\bar{r}))}(\bar{r})
=\Phi_{\chi(\beta_q(\bar{r}))g}(r)
=\Phi_
{\chi(\beta_q(\bar{r}))g\chi(\beta_q(r))^{-1}}(\Phi_{\chi(\beta_q(r))}(r))
=\Phi_h(\Phi_{\chi(\beta_q(r))}(r)),
\]
where $h : = \chi(\beta_q(\bar{r}))g\chi(\beta_q(r))^{-1})$.
Since $\Phi_{\chi(\beta_q(\bar{r}))}(\bar{r}),
\Phi_{\chi(\beta_q(r))}(r)\in S_q$, we have $h\in G_q$ by
$(\rm{S}2)$. So we obtain that
$\Phi_h(\Phi_{\chi(\beta_q(r))}(r))=\Phi_{\chi(\beta_q(r))}(r)$,
by the consequence $(1)$ of the Slice Theorem combined with the
fact that all isotropy groups are conjugated.

Since the value $B_{q_\omega}(r)$ does not depend on the choice of
the representative in the equivalence class of $r$ intersected
with $U$, we can define the map $\mathcal{B}_{q_\omega} :
\mathcal{U} \longrightarrow N^s_{q_\omega}$ such that
$B_{q_\omega}=\mathcal{B}_{q_\omega}\circ\pi$;  it is continuous
by construction.

We now show that $\mathcal{B}_{q_\omega}$ is injective. So let
$\upsilon,\tau\in \mathcal {U}$ such that
$\mathcal{B}_{q_\omega}(\upsilon)=\mathcal{B}_{q_\omega}(\tau)$.
Then choosing any $r_\upsilon\in\pi^{-1}(\upsilon)\cap U$ and
$r_\tau\in\pi^{-1}(\tau)\cap U$ we have:
\begin{align*}
\mathcal{B}_{q_\omega}(\upsilon)=\mathcal{B}_{q_\omega}(\tau)&\Rightarrow
B_{q_\omega}(r_\upsilon)=B_
{q_\omega}(r_\tau)\\
&\Rightarrow
\Phi_{\chi(\beta_{q_\omega}(r_\upsilon))}(r_\upsilon)=\Phi_{\chi
(\beta_{q_\omega}(r_\tau))}(r_\tau)\\
&\Rightarrow \pi(r_\upsilon)=\pi(r_\tau)\\
&\Rightarrow \upsilon=\tau.\;\;\;\;\blacksquare
\end{align*}

\begin{theorem}
\label{TM} Let $\omega\in Q^s/G, s> s_0+1$,  $q_\omega\in
\pi^{-1}(\omega)$, and  $(\psi,V)$ a chart of $Q^s$ at
$q_\omega$. We can choose the neighborhood $\mathcal{U}$
of $\omega$ in $Q^s/G$ such that the map
\[
\mathcal{B}_{q_\omega} : \mathcal{U}\subset Q^s/G \longrightarrow
\mathcal{B}_{q_\omega}(\mathcal{U})\subset N^s_{q_\omega}
\]
is an homeomorphism between open subsets with inverse given by
$\pi\circ\psi^{-1}$. So $Q^s/G$ is a topological manifold (with
the quotient topology) modeled on the Banach space
$N^s_{q_\omega}$.
\end{theorem}
\textbf{Proof.} We proceed in several steps.

\noindent $(1)$ By Lemma \ref{A}, for $V$ sufficiently small, the
continuous map $\pi\circ\psi^{-1} : N_ {q_\omega}^s\cap \psi(V)
\longrightarrow Q^s/G$ is injective.\\
$(2)$ By Lemma \ref{B} there is a neighborhood $\mathcal{U}$ of
$\omega$ such that $\mathcal{B}_ {q_\omega} :
\mathcal{U}\longrightarrow \mathcal{B}_{q_\omega} (\mathcal{U})$
is continuous and bijective. We can choose $\mathcal{U}$ such
that
$\mathcal{B}_{q_\omega} (\mathcal{U})\subset N_ {q_\omega}^s\cap \psi(V)$.\\
$(3)$ So we can compose the maps
\[
\mathcal{U}\stackrel{\mathcal{B}_{q_\omega}} {\hbox to
45pt{\rightarrowfill}}
\mathcal{B}_{q_\omega}(\mathcal{U})\stackrel{\pi\circ\psi^{-1}}
{\hbox to 45pt{\rightarrowfill}}
\pi\left(\psi^{-1}\left(\mathcal{B}_{q_\omega}(\mathcal{U})\right)\right)
\]
to get
\[
\left(\pi\circ\psi^{-1}\circ\mathcal{B}_{q_\omega}\right)(\upsilon)
=\pi\left(\psi^{-1}\left(B_{q_\omega}(r_\upsilon)\right)\right)
=\pi\left(\psi^{-1}\left(\psi\left
(\Phi_{\chi(\beta_{q_\omega}(r_\upsilon))}(r_\upsilon)\right)\right)\right)
=\upsilon,
\]
that is, $\pi\circ\psi^{-1}\circ\mathcal{B}_{q_\omega}$ is the
identity  map on  $\mathcal{U}$. Thus we conclude that:
\begin{itemize}
\item
$\pi\left(\psi^{-1}\left(\mathcal{B}_{q_\omega}(\mathcal{U})\right)\right)
=\mathcal{U}$, so it is an open subset of $Q^s/G$, \item
$\mathcal{B}_{q_\omega}(\mathcal{U})=(\pi\circ\psi^{-1})^{-1}\left(\mathcal{U}
\right)$,
so it is an open subset of $N_ {q_\omega}^s$, and \item
$\pi\circ\psi^{-1}=\left(\mathcal{B}_{q_\omega}\right)^{- 1}$.
\end{itemize}
This shows that $Q^s/G$ is a topological  manifold (with the
quotient topology) modeled on the Banach space $N^s_{q_\omega}$.
$\;\;\;\;\blacksquare$

\medskip

Now a natural question arises: is $Q^s/G$ a differentiable
manifold relative to the naturally induced differentiable
structure? The answer is no. By contradiction, suppose it is the
case, so $\pi : Q^s\longrightarrow Q^s/G$ is a locally trivial
fiber bundle and we can take the tangent map $T_q\pi : T_qQ^s
\longrightarrow T_ {\pi(q)}(Q^s/G)$ whose kernel is the vertical
subspace $V_qQ^s=\operatorname{ker}(T_q\pi)$. We know that the vertical
subspace is generated by the infinitesimal generators at $q$, that
is, $V_qQ^s=\{\xi_{Q^s}(q)\;|\;\xi\in\mathfrak{g}\}$. This is a
contradiction because, in general, $\xi_{Q^s}(q) \notin T_qQ^s$
(recall that we only have $\xi_{Q^s}(q)\in T_qQ^ {s-1}$). Since $Q^s/G$
is not a differentiable manifold, we cannot define the tangent
bundle or $C^1$ curves in the usual way. To overcome this difficulty,
we will use that changes of charts (which are not differentiable in the
usual sense) are not only homeomorphisms, but are also differentiable
relative to a weaker topology  (in fact $H^{s-1}$). This agrees with
the fact that our vertical subspace is a  subset of
$TQ^{s-1}|Q^s$.

For simplicity we suppose that the action $\Phi$ is free, but all
the following results are still true in the case when all the
isotropy groups are conjugated.

In the following definition, we give a generalization of the
notion of $C^1$ curves in $Q^s$ and $N^s_q$. Differentiating
these curves we obtain a generalization of the notion of tangent
vectors to $Q^s$ and $N^s_q$.

\begin{definition}
\label{weakly_diff} Let $s> s_0+1$.
\begin{enumerate}
\item[{\rm (i)}] For $I$ open in $\mathbb{R}$ we define the set
\[
C^1_W(I,Q^s):=C^0(I,Q^s)\cap C^1(I,Q^{s-1})
\]
of continuous curves in $Q^s$ which are continuously
differentiable with respect to the $Q^{s-1}$ differentiable
structure. The set $C^1_W(I,N^s_q)$ is defined in the same way.
Such curves will be called \textbf{weakly-differentiable}.
\item[{\rm (ii)}] For $q\in Q^s$, the \textbf{weak tangent space}
of $Q^s$ at $q$ is defined by
\[
T^W_qQ^s:=\{\dot{d}(0)\in T_qQ^{s-1}\,|\,d\in C^1_W(I,Q^s),
d(0)=q\}.
\]
The set $T^W_{v_q}N^s_q$ is defined in the same way.
\item[{\rm(iii)}] The \textbf{weak tangent bundle} of $Q^s$ is
defined by
\[
T^WQ^s:=\bigcup_{q\in Q^s}T^W_qQ^s.
\]
The set  $T^WN^s_q$ is defined in the same way.
\end{enumerate}

\end{definition}

Note that $T^W_qQ^s$ is a vector space and that we have the
inclusions $T_qQ^s\subset T^W_qQ^s\subset T_qQ^{s-1}$, for all
$q\in Q^s$. The same is true for $N^s_q$. Since $N^s_q$ is a
vector space, $T^W_{v_q}N^s_q$ does not depend on $v_q$ and we
prefer the notation $(N^s_q)^W:=T^W_{v_q}N^s_q$. The previous
inclusions become in this case $N^s_q \subset (N^s_q)^W \subset
N^{s-1}_q$ and the tangent bundle is $T^WN^s_q=N^s_q\times
(N^s_q)^W$. We endow $T^W_qQ^s$ and $(N^s_q)^W$ with the topology
of $T_qQ^{s-1}$ and $N^{s-1}_q$, respectively. Note that for all
$q\in Q^s$ and $\xi\in\mathfrak{g}$, we have $\xi_{Q^s}(q)\in
T^W_qQ^s$, since
$\xi_{Q^s}(q)=\left.\frac{d}{dt}\right|_{t=0}\Phi_{\operatorname{exp}(t\xi)}(q)
$
and the curve $t\mapsto\Phi_{\operatorname{exp}(t\xi)}(q)$ is in
$C^1_W(I,Q^s)$ by the
hypothesis $\ref{phi_c_one}$ on the $G$-action.\\

Now we consider the regularity of the map $B_q$. To be more
precise we will use the notation $B^s_q$ instead of $B_q$. Recall
that $B^s_q$ is given by
\[
B^s_q : U\subset Q^s\longrightarrow N^s_q,\;\;B^
s_q(r):=\varphi^s(\Phi_{\beta_q(r)}(r)).
\]
Note that $\beta_{q}$ is $C^1$ by Theorem \ref{IFThm} and that
$\Phi$ is $C^1$ as a map with values in $Q^{s-1}$ by the working
hypothesis on the action \eqref{phi_c_one}. Hence $r \in U \mapsto
\Phi_{\beta_q(r)}(r) \in Q^{s-1}$ is a $C ^1$ map. Since we can
think of the chart $\varphi^s$ of $Q^s$ at $q$ as the
restriction of a chart $\varphi^{s-1}$ of $Q^{s-1}$ at $q$, we can
choose $U$ such that
\[
B^s_q\in C^1\left(U,N^{s-1}_q\right).
\]
Interpreted this way, we can compute the tangent map at $r$ in the
direction $v_r\in TQ^s|U$. We have by Lemma \ref{B}
\begin{equation}
\label{tangentofB} T_rB^s_q(v_r) =T_{\Phi_{\beta_q(r)}(r)}
\varphi^{s-1}\left(T_r\Phi_{\beta_{q}(r)}
(v_r)+T_{\beta_q(r)}\Phi^r(T_r\beta_{q}(v_r))\right).
\end{equation}
Now we will show that this tangent map
\[
TB^s_q : TQ^s|U \longrightarrow N^s_q\times N^{s-1}_q
\]
makes sense and is continuous on the larger space $TQ^{s-1}|U$.
Indeed, by the assumption (\ref{phi_g_smooth}), the term
$T_r\Phi_{\beta_{q}(r)} (v_r)$ makes sense for all $v_r\in
TQ^{s-1}|U$. It remains to show that $T_r\beta_{q}(v_r)$ makes
sense for $v_r\in TQ^{s-1}|U$.

By Theorem \ref{IFThm} we have $\mathcal{S}_q(\beta_q(r),r)=0$ for
$r$ in a neighborhood of $q$. Differentiating this expression relative
to $r $, we obtain
\[
T_{\beta_q(r)}[\mathcal{S}_q(\_,r)](T_r\beta_q(v_r))+T_r[\mathcal{S}_q(\beta_q
(r),\_)](v_r)=0.
\]
By Theorem \ref{IFThm} again, we know that
$T_{\beta_q(r)}[\mathcal{S}_q(\_,r)]$ is invertible in a
neighborhood of $q$, so we obtain
\[
T_r\beta_q(v_r)=-\Big{(}T_{\beta_q(r)}[\mathcal{S}_q(\_,r)]\Big{)}^{-1}T_r
[\mathcal{S}_q(\beta_q(r),\_)](v_r).
\]
Thus it suffices to show that
$T_r[\mathcal{S}_q(\beta_q(r),\_)](v_r)$ is well-defined for
$v_r\in TQ^{s-1}|U$. This is true because
\eqref{s_q_useful_formula} implies the formula
\[
T_r[\mathcal{S}_q(g,\_)](v_r)=\sum_{i=1}^n\gamma(q)\left(T_{\Phi_g(r)}\varphi^
{s-1}(T\Phi_g(v_r)),E_i(q)\right)e_i.
\]
So the tangent map $TB^s_q$ makes sense on $TQ^{s-1}|U$ and
therefore on $T^WQ^s|U$. We will denote by
\[
T^WB^s_q : T^WQ^s|U \longrightarrow N^s_q\times N^{s-1}_q
\]
this extension. Remark that if $s> s_0+2$ we have
$T^WB^s_q=TB^{s-1}_q$ on $T^WQ^s|U$. The following Lemma shows
that $T^WB^s_q$ takes values in $T^WN^s_q=N^s_q\times
(N_q^s)^W\subset N^s_q\times N^{s-1}_q$, that is,
\[
T^WB^s_q : T^WQ^s|U \longrightarrow T^WN^s_q=N^s_q\times
(N^s_q)^W.
\]

\begin{lemma}
\label{chain rule} Let $q\in Q^s, s> s_0+2$, and $U$ neighborhood
of $q$ in $Q^s$ such that $B^s_q$ is defined on $U$. Let $d\in
C^1_W(I,Q^s)$ such that $d(t)\in U$ for all $t\in I$. Then
$B^s_q\circ d\in C^1_W(I,N^s_q)$ and we have
\begin{equation}
\label{chain rule} \frac{d}{dt}(B^s_q\circ
d)(t)=T^W_{d(t)}B^s_q(\dot{d}(t)).
\end{equation}
\end{lemma}
\textbf{Proof.} It is clear that $B^s_q\circ d\in C^0(I,N^s_q)$ so
we have to show that $B^s_q\circ d\in C^1(I,N^{s-1}_q)$. Recall
that $B_q(d(t))=\varphi^s(\Phi_{\beta_q(d(t))}(d(t)))$. For all
$u\in I$, by \eqref{phi_g_smooth}, we have
\[
\begin{array}{ccccc}
I&\stackrel{C^1}{\longrightarrow}&Q^{s-1}&\stackrel{C^\infty}
{\longrightarrow} &Q^
{s-1}\\
t&\longmapsto                    & d(t)    &\longmapsto

&\Phi_{\beta_q(d(u))}(d(t)).
\end{array}
\]
Since  $d \in C ^1(I, Q^{s-1})$, using Theorem \ref{IFThm} (for
$s-1$, which is allowed since we assume $s> s_0 +2$), and
\eqref{phi_c_one}, we get for all $t\in I$
\[
\begin{array}{ccccccc}
I&\stackrel{C^1}{\longrightarrow}&Q^{s-1}&\stackrel{C^1}{\longrightarrow}

&G          &\stackrel{C^1}{\longrightarrow} &Q^{s-1}\\
u&\longmapsto                    & d(u)    &\longmapsto

&\beta_q(d(u))&\longmapsto &\Phi_{\beta_q(d(u))}(d(t))
\end{array}
\]
So $B^s_q\circ d\in C^1(I,N^{s-1}_q)$.

To  prove formula $\eqref{chain rule}$ recall that  $d\in
C^1(I,Q^{s-1})$ and $B^{s-1}_q\in C^1(U,N^{s-2}_q)$, so by the
chain rule we have $B^{s-1}_q\circ d\in C^1(I,N^{s-2}_q)$ and
\[
\frac{d}{dt}(B^{s-1}_q\circ d)(t)=T_{d(t)}B^{s-1}_q(\dot{d}(t)).
\]
Since $d(t)\in U$ we have $B^{s-1}_q\circ d=B^s_q\circ d$ and, as
was discussed before,
$T_{d(t)}B^{s-1}_q(\dot{d}(t))=T_{d(t)}^WB^s_q(\dot{d}(t))$. So we
conclude that
\[
\frac{d}{dt}(B^s_q\circ d)(t)=T^WB^s_q(\dot{d}(t))
\]
where $B^s_q\circ d$ is derived as a curve in $C^1(I,N^{s-2}_q)$.
But we know that we have in fact $B^s_q\circ d\in
C^1(I,N^{s-1}_q)$ so in the preceding formula $B^s_q\circ d$ is
derived as a curve in $C^1(I,N^{s- 1}_q).\;\;\;\;\blacksquare$

\medskip

Now we state an important result of this section. It says in which
sense the change of charts is differentiable. This will be useful
in the construction of a tangent bundle of $Q^s/G$.

\medskip

\begin{theorem}
\label{change of charts} Let $\omega,\upsilon\in Q^s/G, s>s_0+2$,
and $(\mathcal{U}_1,\mathcal{B}_{q_\omega}),
(\mathcal{U}_2,\mathcal{B}_{q_\upsilon})$ be two charts of $Q^s/G$
such that $\mathcal{U}_1\cap\mathcal{U}_2\neq \varnothing$. Then
the change of charts
\[
F:=\mathcal{B}_{q_\upsilon}\circ \mathcal{B}_{q_\omega}^{-1} :
\mathcal{B}_{q_\omega}(\mathcal{U}_1\cap\mathcal{U}_2)\subset
N^s_{q_\omega}\longrightarrow
\mathcal{B}_{q_\upsilon}(\mathcal{U}_1\cap\mathcal{U}_2)\subset
N^s_{q_\upsilon}
\]
is $C^1$ as a map with values in $N^{s-1}_{q_\upsilon}$.
\end{theorem}
\textbf{Proof.} By construction of the charts, see Lemma \ref{B},
we have $\mathcal{U}_1=\pi(U_1)$ and $\mathcal{U}_1=\pi(U_1)$
where $U_1$ is a neighborhood of $q_\omega$ and $U_2$ a
neighborhood of $q_\upsilon$. Recall that
$\mathcal{B}_{q_\upsilon}$ is defined such that
$\mathcal{B}_{q_\upsilon}\circ\pi=B_{q_\upsilon} :
U_2\longrightarrow N^s_{q_\upsilon}$, where
$B_{q_\upsilon}(r)=\psi^s(\Phi_{\beta_{q_\upsilon}(r)}(r))$.

Let $v_{q_\omega}\in
\mathcal{B}_{q_\omega}(\mathcal{U}_1\cap\mathcal{U}_2)$ and let
$\tau:=\mathcal{B}_{q_\omega}^{-1}(v_{q_\omega})=(\pi\circ(\varphi^s)^{-1})(v_
{q_\omega})\in\mathcal{U}_1\cap\mathcal{U}_2$.
Since $\pi((\varphi^s)^{-1}(v_{q_\omega}))\in
\mathcal{U}_2=\pi(U_2)$, there exists $g_0\in G$ such that
$\Phi_{g_0}((\varphi^s)^{-1}(v_{q_\omega}))\in U_2$. By continuity
of $\Phi_{g_0}$, there exists a  neighborhood $V$ of
$(\varphi^s)^{-1}(v_{q_\omega})$ such that $\Phi_{g_0}(V)\in U_2$.
Since $v_{q_\omega}\in
\mathcal{B}_{q_\omega}(\mathcal{U}_1\cap\mathcal{U}_2)$ we can
assume that $\varphi^s(V)\cap N^s_{q_\omega}\subset
\mathcal{B}_{q_\omega}(\mathcal{U}_1\cap\mathcal{U}_2)$.

So for all $u_{q_\omega}$ in the neighborhood $\varphi^s(V)\cap
N^s_{q_\omega}$ of $v_{q_\omega}$, we have:
\begin{align*}
F(u_{q_\omega})&=\mathcal{B}_{q_\upsilon}\left(\mathcal{B}_{q_\omega}^{-1}
(u_{q_\omega})\right)\\
&=\mathcal{B}_{q_\upsilon}\left(\pi((\varphi^s)^{-1}(u_{q_\omega}))\right)\\
&=\mathcal{B}_{q_\upsilon}\left(\pi(\Phi_{g_0}((\varphi^s)^{-1}(u_{q_\omega})))
\right)\\
&=B_{q_\upsilon}(\Phi_{g_0}((\varphi^s)^{-1}(u_{q_\omega}))),
\text{ since $\Phi_{g_0}((\varphi^s)^{-1}(u_{q_\omega}))\in
U_2$}\\
&=\varphi^s\left(\Phi_{\beta_{q_\upsilon}(\Phi_{g_0}((\varphi^s)^{-1}(u_
{q_\omega})))}\left(\Phi_{g_0}((\varphi^s)^{-1}(u_{q_\omega}))\right)\right).
\end{align*}
An inspection of this formula (using the assumption
(\ref{phi_c_one}) and the fact that the chart $\varphi^s$ of
$Q^s$ can be seen as the restriction of a chart $\varphi^{s-1}$
of
$Q^{s-1}$) shows that, in a neighborhood of $v_{q_\omega}$, $F$ is
$C^1$ as a map with values in $N^{s-1}_{q_\upsilon}$. Doing that
for all
$v_{q_\omega}\in\mathcal{B}_{q_\omega}(\mathcal{U}_1\cap\mathcal{U}_2)$
we obtain the desired result.$\;\;\;\;\blacksquare$

\medskip

Now we consider the tangent map to $F$:
\[
TF : \mathcal{B}_{q_\omega}(\mathcal{U}_1\cap\mathcal{U}_2)\times
N^s_{q_\omega} \rightarrow
\mathcal{B}_{q_\upsilon}(\mathcal{U}_1\cap\mathcal{U}_2)\times
N^{s-1}_{q_\upsilon},\;\;TF(u_
{q_\omega},v_{q_\omega}):=\left(F(u_{q_\omega}),D\,F(u_{q_\omega})(v_
{q_\omega})\right).
\]
As it was the case for $B_q$, $TF$ makes sense and is continuous
on the larger space
$\mathcal{B}_{q_\omega}(\mathcal{U}_1\cap\mathcal{U}_2)\times
(N^s_{q_\omega})^W$. We will denote by
\[
T^WF :
\mathcal{B}_{q_\omega}(\mathcal{U}_1\cap\mathcal{U}_2)\times
(N^s_{q_\omega})^W \longrightarrow
\mathcal{B}_{q_\upsilon}(\mathcal{U}_1\cap\mathcal{U}_2)\times
N^{s-1}_{q_\upsilon}
\]
this extension. As in Lemma \ref{chain rule} we can show that for
$d\in C^1_W (I,N^s_{q_x})$, $s> s_0+2$, such that $d(t)\in U$ for
all $t\in I$, we have $F\circ d\in C^1_W(I,N^s_{q_y})$ and
\begin{equation}
\label{chain rule for F} \frac{d}{dt}(F\circ
d)(t)=T^W_{d(t)}F(\dot{d}(t)).
\end{equation}
Thus $T^WF$ takes values in
$\mathcal{B}_{q_\upsilon}(\mathcal{U}_1\cap\mathcal{U}_2)\times
(N^s_{q_\upsilon})^W$, that is
\[
T^WF :
\mathcal{B}_{q_\omega}(\mathcal{U}_1\cap\mathcal{U}_2)\times
(N^s_{q_\omega})^W \longrightarrow
\mathcal{B}_{q_\upsilon}(\mathcal{U}_1\cap\mathcal{U}_2)\times
(N^s_{q_\upsilon})^W.
\]
Of course Theorem \ref{change of charts} and all previous remarks
are valid for
$F^{-1}=\mathcal{B}_{q_\omega}\circ(\mathcal{B}_{q_\upsilon})^{-1}
: \mathcal{B}_{q_\upsilon}(\mathcal{U}_1\cap\mathcal{U}_2)\subset
N^s_ {q_\upsilon} \longrightarrow
\mathcal{B}_{q_\omega}(\mathcal{U}_1\cap\mathcal{U}_2)\subset
N^s_{q_\omega}$ and we can construct the tangent map
$T^W(F^{-1})$. We have naturally the following property.

\begin{theorem} Under the hypothesis of the preceding theorem,
\[
T^WF :
\mathcal{B}_{q_\omega}(\mathcal{U}_1\cap\mathcal{U}_2)\times
(N^s_{q_\omega})^W \longrightarrow
\mathcal{B}_{q_\upsilon}(\mathcal{U}_1\cap\mathcal{U}_2)\times
(N^s_{q_\upsilon})^W
\]
is an homeomorphism whose inverse is given by $T^W(F^{-1})$.
\end{theorem}
\textbf{Proof.} As we remarked before, $T^WF$ is clearly
continuous on
$\mathcal{B}_{q_\omega}(\mathcal{U}_1\cap\mathcal{U}_2)\times
(N^s_{q_\omega})^W$. Let $d\in C^1_W(I,N^s_{q_x})$, so we have
$F\circ d\in C^1_W(I,N^s_{q_y})$ and we can compute:
\begin{align*}\dot{d}(t)&=\frac{d}{dt}(F^{-1}\circ F\circ d)(t)\\
&=T^W_{F(d(t))}(F^{-1})\left(\frac{d}{dt}(F\circ
d)(t)\right)\;\;\text{ with formula (\ref
{chain rule for F}), still valid for $F^{-1}$}\\
&=T^W_{F(d(t))}(F^{-1})\left(T^W_{d(t)}F(\dot{d}(t))\right)\;\;\text{
with formula (\ref{chain rule for F}) again}.\;\;\;\;\blacksquare
\end{align*}

\begin{definition}
\label{weakly differentialble curves}
For $s> s_0+1$ and $I$ open in $\mathbb{R}$ we define
the set $C^1_W(I,Q^s/G)$ of \textbf{weakly-differentiable} curves
in $Q^s/G$ as follows: \quote{$c\in C^1_W(I,Q^s/G)$ if and only if
for all $t_0\in I$ there exists $\varepsilon >0$ and $d\in
C^1_W(]t_0-\varepsilon,t_0+\varepsilon[,Q^s)$ such that}
\[
c=\pi\circ d\;\;\text{ on }\;\;]t_0-\varepsilon,t_0+\varepsilon[ \subset
I.
\]
\end{definition}

The following theorem shows that in a chart the notion of
weakly-differentiable curves in $Q^s/G$ coincides with the notion
of weakly-differentiable curves in $N^s_q$ given in Definition
\ref{weakly_diff}.

\begin{theorem}
\label{equivalent definition} Let $s> s_0+2$ and $I$ be open in
$\mathbb{R}$. Then $c\in C^1_W(I,Q^s/G)$ if and only if for all
charts we have $\mathcal{B}_q\circ c\in C^1_W(I',N^s_q)$, where
$I'$ is an open subset of $I$ such that $\mathcal{B}_q\circ c$ is
well-defined.
\end{theorem}
\textbf{Proof.} If $c\in C^1_W(I,Q^s/G)$ then, by definition,
$c=\pi\circ d$ on $]t_0-\varepsilon,t_0+\varepsilon[$, where $d\in
C^1_W(]t_0-\varepsilon,t_0+\varepsilon[,Q^s)$. So we obtain that
$\mathcal{B}_q\circ c=B_q\circ d\in C^1_W(I',N^s_q)$ by Lemma
\ref{chain rule}.

Conversely, let $t_0\in I$, $\omega:=c(t_0)$ and $q\in
\pi^{-1}(\omega)$. We have $\mathcal{B}_q\circ c\in
C^1_W(]t_0-\varepsilon,t_0+\varepsilon[,N^s_q) $. Let
$d:=(\varphi^s)^{-1}\circ \mathcal{B}_q\circ c$. Then $d\in C^1_W
(]t_0-\varepsilon,t_0+\varepsilon[,Q^s)$ and $\pi\circ
d=\pi\circ(\varphi^s)^{-1}\circ\mathcal{B}_q\circ c=c$. Doing that
for each $t_0\in I$ shows that $c\in
C^1_W(I,Q^s/G).\;\;\;\;\blacksquare$

\medskip
Now we are ready to construct the tangent space of the isotropy
strata. This will be done using the notion of
weakly-differentiable tangent
 curves at one point.

\begin{definition} Let $\omega\in Q^s/G, s> s_0+2$, and $c_1,c_2\in C^1_W
(I,Q^s/G)$ such that $c_1(0)=c_2(0)=\omega$. We define the
following relation:
\[
c_1\stackrel{\omega}{\sim}c_2\Longleftrightarrow
\left.\frac{d}{dt}\right|_{t=0} (\mathcal{B}_q\circ
c_1)(t)=\left.\frac{d}{dt}\right|_{t=0}(\mathcal{B}_q\circ c_2)(t)
\]
where $\mathcal{B}_q$ is any chart at $\omega$.
\end{definition}

We emphasize that $\mathcal{B}_q\circ c_i$ is a curve in $N^s_q$
but it is derived as a $C^1$ curve in $N^{s-1}_q$ by Theorem \ref
{equivalent definition}.

\begin{lemma}
\label{equivB} In the preceding definition, the relation
$\stackrel{\omega}{\sim}$ does not depend on the choice of the
chart $\mathcal{B}_q$. Thus $\stackrel{\omega}{\sim}$ is an
equivalence relation on the set $C^1_{W,\omega}(I,Q^s/G):=\{c\in
C^1_W(I,Q^s/G)\;|\;c(0)=\omega\}, s>s_0+2$.
\end{lemma}
\textbf{Proof.} The proof is like the standard one but we shall do
it because our curves are differentiable in a weaker sense. Let's
suppose that $\left.\frac{d}{dt}\right|_{t=0}(\mathcal{B}_q\circ
c_1)(t)=\left.\frac{d}{dt}\right|_{t=0}(\mathcal{B}_q\circ
c_2)(t)$ and let $\mathcal{B}_r$ be another chart at $\omega$. We
have
\begin{align*}
\left.\frac{d}{dt}\right|_{t=0}(\mathcal{B}_q\circ
c_1)(t)&=\left.\frac{d}
{dt}\right|_{t=0}(\mathcal{B}_q\circ(\mathcal{B}_r)^{-1}\circ\mathcal{B}
_r\circ c_1)(t)\\
&=\left.\frac{d}{dt}\right|_{t=0}(F\circ (\mathcal{B}_r\circ c_1))(t)\\
&=T^W_{\mathcal{B}_r(\omega)}F\left(\left.\frac{d}{dt}\right|_{t=0}(\mathcal{B}
_r\circ
c_1)(t)\right) \;\;\text{ because of formula }\eqref{chain rule
for F}.
\end{align*}
Doing the same for $c_2$ and using the bijectivity of $T^WF$ we
obtain that:
\[
\left.\frac{d}{dt}\right|_{t=0}(\mathcal{B}_r\circ
c_1)(t)=\left.\frac{d}{dt} \right|_{t=0}(\mathcal{B}_r\circ
c_2)(t).\;\;\;\;\blacksquare
\]

\begin{definition}\label{definition_of_weak} Let $s> s_0+2$.
\begin{enumerate}
\item[{\rm (i)}] Let $\omega\in Q^s/G$. The \textbf{weak tangent
space} of $Q^s/G$ at $\omega$ is defined by
\[
T_\omega^W(Q^s/G):=C^1_{W,\omega}(I,Q^s/G)/\stackrel{\omega}{\sim}.
\]
We will denote by $v_\omega$, $[c]_\omega$, or $\dot{c}(0)$ the
elements of $T_\omega^W(Q^s/G)$. \item[{\rm (ii)}] The
\textbf{weak tangent bundle} of $Q^s/G$ is defined by
\[
T^W(Q^s/G):=\bigcup_{\omega\in Q^s/G}T_\omega^W(Q^s/G).
\]
\item[{\rm (iii)}] Let $(\mathcal{U},\mathcal{B}_{q_\omega})$ be a
chart of $Q^s/G$. The \textbf{weak tangent map} of
$\mathcal{B}_{q_\omega}$ is defined by
\[
T^W\mathcal{B}_{q_\omega} : T^W(Q^s/G)|\mathcal{U}\longrightarrow
\mathcal{B}_{q_\omega}(\mathcal{U})\times (N^s_{q_\omega})^W,
\]
\[
T^W_\upsilon\mathcal{B}_{q_\omega}([c]_\upsilon):=\left.\frac{d}{dt}\right|_
{t=0}(\mathcal{B}_{q_\omega}\circ
c)(t).
\]
\item[{\rm (iv)}] Let $(\mathcal{U},\mathcal{B}_{q_\omega})$ be a
chart of $Q^s/G$. The \textbf{weak tangent map} of
$(\mathcal{B}_{q_\omega})^{-1}$ is defined by
\[
T^W(\mathcal{B}_{q_\omega})^{-1} :
\mathcal{B}_{q_\omega}(\mathcal{U})\times
(N^s_{q_\omega})^W\longrightarrow T^W(Q^s/G)|\mathcal{U},
\]
\[
T^W_{v_{q_\omega}}(\mathcal{B}_{q_\omega})^{-1}(v_{q_\omega},w_{q_\omega}):=
[(\mathcal{B}_
{q_\omega})^{-1} \circ
d]_{(\mathcal{B}_{q_\omega})^{-1}(v_{q_\omega})}
\]
where $d\in C^1_W(I,N^s_{q_\omega})$ is such that
$\dot{d}(0)=(v_{q_\omega},w_{q_\omega}) $. \item[{\rm (v)}] The
\textbf{weak tangent map} of $\pi$ is defined by
\[
T^W\pi : T^WQ^s \longrightarrow T^W
(Q^s/G),\;\;T^W_r\pi(v_r):=[\pi\circ d]_{\pi(r)}
\]
where $d\in C^1_W(I,Q^s)$ is such that $\dot{d}(0)=v_r$.
\end{enumerate}
\end{definition}

The properties of these weak tangent maps are summarized in the
following statement.

\begin{lemma}
\label{properties_B} We have the following.
\begin{enumerate}
\item[{\rm (i)}] $T^W(\mathcal{B}_{q_\omega})^{-1}$ is well
defined, that is, it does not depend on the choice of the curve
$d$. Moreover:
\[
T^W\mathcal{B}_{q_\omega}\circ T^W(\mathcal{B}_{q_\omega})^{-1}=id
\]
\item[{\rm (ii)}] $T^W\pi$ is well-defined, that is, it does not
depend on the choice of the curve $d$. Moreover for all charts
$\mathcal{B}_{q_\omega}$ we have on $T^WQ^s|U$:
\[
T^W\mathcal{B}_{q_\omega}\circ
T^W\pi=T^W(\mathcal{B}_{q_\omega}\circ \pi) =T^WB_{q_\omega}.
\]
\item[{\rm (iii)}] Under the assumptions of Theorem \ref{change of
charts} we have:
\[
T^WF=T^W\mathcal{B}_{q_\upsilon}\circ
T^W(\mathcal{B}_{q_\omega})^{-1}.
\]
\end{enumerate}
\end{lemma}

\textbf{Proof.} (i) Let $d(t)$ be a curve in $C^1_W(I,N_{q_x}^s)$
such that $\dot{d}(0)=(v_{q_\omega},w_{q_\omega})\in
\mathcal{B}_{q_\omega}(\mathcal{U})\times (N^s_{q_\omega})^W$.
Then we get
\begin{align*}
\left(T^W\mathcal{B}_{q_\omega}\circ
T^W(\mathcal{B}_{q_\omega})^{-1}\right) (\dot{d}(0))  &=
T^W_{(\mathcal{B}
_{q_\omega})^{-1}(v_{q_\omega})}\mathcal{B}_{q_\omega}\left([(\mathcal{B}_
{q_\omega})^{-1}
\circ d]_{(\mathcal{B}_{q_\omega})^{-1}(v_{q_\omega})}\right)\\
&=\left.\frac{d}{dt}\right|_{t=0}\left(\mathcal
{B}_{q_\omega}((\mathcal{B}_{q_\omega})^{-1}(d(t)))\right)=\dot{d}(0)
\end{align*}
by the definitions of $T^W\mathcal{B}_{q_\omega}$ and
$T^W(\mathcal{B}_{q_\omega})^{-1}$. This shows that
$T^W\mathcal{B}_{q_\omega}\circ
T^W(\mathcal{B}_{q_\omega})^{-1}=id$. On the other hand, since
$T^W_{v_{q_\omega}}(\mathcal{B}_{q_\omega})^{-
1}(w_{q_\omega}):=[(\mathcal{B}_{q_\omega})^{-1}\circ
d]_{(\mathcal{B}_{q_\omega}) ^{-1}(v_{q_\omega})} = \left(T^W
\mathcal{B}_{q_\omega} \right)^{-1} (\dot{d}(0))$ we see that this
expression depends only on
$\dot{d}(0)=(v_{q_\omega},w_{q_\omega})$ and not on the curve
$d(t)$ itself.

(ii) For a curve $d(t)$ in $C^1_W(I,N_{q_\omega}^s)$ such that
$d(0)=r\in U$ we have
\[
T^W_{\pi(r)}\mathcal{B}_{q_\omega}\left([\pi\circ
d]_{\pi(r)}\right)=\left.\frac{d}{dt}
\right|_{t=0}\left(\mathcal{B}_{q_\omega}(\pi(d(t)))\right)=\left.\frac{d}{dt}
\right|_{t=0}\left(B_{q_\omega}(d(t))\right)=T^W_rB_{q_\omega}(\dot{d}(0))
\]
by the definition of $T^W\mathcal{B}_{q_\omega}$ and using Lemma
\ref{chain rule}. Since the right hand side does not depend on
$d(t)$, it follows that $T^W\pi(v_r):=[\pi\circ d]_{\pi(r)}$ also
does not depend on $d(t)$ satisfying $\dot{d}(0)=v_r$. In addition,
the formula shows that $T^W\mathcal{B}_{q_\omega}\circ T^W\pi=T^W
(\mathcal{B}_{q_\omega}\circ
\pi)=T^WB_{q_\omega}.$\\
(iii) This follows directly from the definitions.
$\;\;\;\;\blacksquare$

\medskip

Note that we can give $T_\omega^W(Q^s/G)$ the structure of a
vector space as it is done in the standard case.

\begin{theorem} Let $s> s_0+2$. Then:
\begin{enumerate}
\item[{\rm (i)}] $T^W(Q^s/G)$ is a topological manifold modeled on
the space $N^s_q\times (N^s_q)^W$ with $q\in Q^s$. \item[{\rm
(ii)}] $T^W\mathcal{B}_{q_\omega},
T^W(\mathcal{B}_{q_\omega})^{-1}$, and $T^W\pi$ are continuous
maps.
\end{enumerate}
\end{theorem}
\textbf{Proof.} (i) By the third part of the previous lemma,
$(T^W\mathcal{U},T^W\mathcal{B}_q)$ are topological charts of $T^W
(Q^s/G)$, where $T^W\mathcal{U}:=T^W(Q^s/G)|\mathcal{U}$. Hence
$T^W(Q^s/G)$ is a topological manifold whose model is the space
$N^s_q\times (N^s_q)^W$ with $q\in Q^s$. (ii) This is
obvious.\;\;\;\;$\blacksquare$

\begin{definition} Let $q\in Q^s, s> s_0+2$.
\begin{enumerate}
\item[{\rm (i)}] The \textbf{vertical subspace} of $T^W_qQ^s$ is
defined by
\[
V_q^WQ^s:=\operatorname{ker}\left(T_q^W\pi\right)
\]
\item[{\rm (ii)}] The \textbf{horizontal subspace} of $T^W_qQ^s$
relative to the metric $\gamma$ is defined by
\[
H_q^WQ^s:=\left(V^W_qQ^s\right)^{\perp}
\]
where $\perp$ means the othogonal complement with respect to the
inner product $\gamma(q)$.
\end{enumerate}
\end{definition}

The following lemma shows that our construction of tangent bundles
in the situation of a non-smooth action is a good generalization
of the standard situation.

\begin{lemma} Let $q\in Q^s, s> s_0+2$. Then we have:
\begin{enumerate}
\item[{\rm (i)}] $V_q^WQ^s=\{\xi_{Q^s}(q)
\;|\;\xi\in\mathfrak{g}\}$. \item[{\rm (ii)}]
$H^W_qQ^s=(N^s_q)^W=\{v_q\in T_q^WQ^s\;|\;\gamma(q)
(v_q,\xi_{Q^s}(q))=0, \forall\,\xi\in\mathfrak{g}\}$. \item[{\rm
(iii)}] The projections associated to the decomposition
$T^W_qQ^s=V_q^WQ^s\oplus H^W_qQ^s$ are given by
\[
\operatorname{ver}_q : T^W_qQ^s \longrightarrow
V^W_qQ^s,\;\;\operatorname
{ver}_q(v_q)=\sum_{i=1}^n\gamma(q)\left(v_q,E_i(q)\right)E_i(q)
\]
\[
\operatorname{hor}_q : T^W_qQ^s \longrightarrow
H^W_qQ^s,\;\;\operatorname
{hor}_q(v_q):=v_q-\operatorname{ver}_q(v_q),
\]
where the basis $(e_1,...,e_n)$ of $\mathfrak{g}$ is chosen such
that $\gamma(q)\left(E_i(q),E_j(q)\right)=\delta_{ij}$. \item[{\rm
(iv)}] Let $\omega\in Q^s/G$ and $q_\omega\in\pi^{-1}(\omega)$.
Then the map
\[
T^W_{q_\omega}\pi : H^W_{q_\omega}Q^s \longrightarrow
T^W_\omega(Q^s/G)
\]
is a continuous linear bijection. Its inverse, called the
\textbf{horizontal-lift}, is denoted by
\[
\operatorname{Hor}_{q_\omega} : T^W_\omega(Q^s/G) \longrightarrow
H^W_{q_\omega}Q^s.
\]
\end{enumerate}
\end{lemma}
\textbf{Proof.} (i) It suffices to prove that
$V_q^WQ^s=\{\sum_{i=1}^n\lambda^iE_i(q)\;|\;\lambda^i\in\mathbb{R}\}$.
Using the chart $\mathcal{B}_q$, \eqref{tangentofB}, $\beta_q
(q)=e$, and $T_q\varphi^s=id$, we have the following equivalences:
\begin{align*}
v_q\in V^W_qQ^s&\Longleftrightarrow
T^W_q\pi(v_q)=0_{\pi(q)}\Longleftrightarrow T^W_{\pi(q)}\mathcal{B}_q(T^W_q\pi
(v_q))=(\mathcal{B}_q(\pi(q)),0_q)\\
&\Longleftrightarrow
T^W_qB_q(v_q)=(\mathcal{B}_q(\pi(q)),0_q)\Longleftrightarrow T_q\Phi_e(v_q)
+T_e\Phi^q(T_q\beta_q(v_q)))=0_q\\
&\Longleftrightarrow
v_q+T_e\Phi^q(T_q\beta_q(v_q))=0_q\Longleftrightarrow
v_q+T_e\Phi^q\left(\sum_{i=1}^nT_q\beta_q(v_q)^ie_i\right)=0_q,
\end{align*}
where $T_q\beta_q(v_q)^i$ are the components of $T_q\beta_q(v_q)$
relative to the basis $(e_1,...,e_n)$. Using that
\[
T_e\Phi^q\left(\sum_{i=1}^nT_q\beta_q(v_q)^ie_i\right)=\sum_{i=1}^nT_q\beta_q
(v_q)^iT_e\Phi^q(e_i)=\sum_{i=1}^nT_q\beta_q(v_q)^iE_i(q)
\]
we obtain that $v_q=\sum_{i=1}^n\lambda^iE_i(q)$.

Conversely, for $v_q=\sum_{i=1}^n\lambda^iE_i(q)$, we have
$T^W_q\pi(v_q)=\sum_{i=1}^n\lambda^iT_q^W\pi(E_i(q))$ and
\[
T_q^W\pi(E_i(q))=T_q^W\pi\left(\left.\frac{d}{dt}
\right|_{t=0}\Phi_{\operatorname{exp}(te_i)}(q)\right)=\left.\frac{d}{dt}
\right|_{t=0}\left(\pi\circ\Phi_{\operatorname{exp}(te_i)}\right)(q)
=\left.\frac{d}{dt}
\right|_{t=0}\pi(q)=0.
\]

(ii) By the definition and (i) we have
\[
H_q^WQ^s:=\left(V^W_qQ^s\right)^{\perp}=\{v_q\in
T_q^WQ^s\;|\;\gamma(q) (v_q,\xi_{Q^s}(q))=0,
\forall\,\xi\in\mathfrak{g}\}.
\]
It remains to show that $H^W_qQ^s=(N^s_q)^W$. The map
\[
T^WB_q : T^WQ^s|U \longrightarrow B_q (U) \times (N^s_q)^W
\subset N^s_q\times (N^s_q)^W
\]
is surjective by Lemma \ref{properties_B} (ii). Using formula
\eqref{tangentofB} we find that $T^WB_q(v_q)=v_q$ for all
$v_q\in H_q^WQ^s$, so we conclude that $v_q\in (N^s_q)^W$.

The other points are obvious.\;\;\;\;$\blacksquare$

\medskip

When the metric $\gamma$ is $G$-invariant, we can define a metric
$\widetilde{\gamma}$ on $Q^s/G$:
\[
\widetilde{\gamma}(\omega)(u_\omega,v_\omega):=\gamma(q_\omega)(\operatorname
{Hor}_{q_\omega}(u_\omega),\operatorname{Hor}_
{q_\omega}(v_\omega))
\]
where $\omega\in Q^s/G$, $u_\omega,v_\omega\in T^W_\omega(Q^s/G)$
and $q_\omega$ is any element in $\pi^{-1}(\omega)$. One can show
that this expression does not depend on the choice of $q_\omega$
in $\pi^{-1}(\omega)$.

\medskip

We now can show in which sense the $(H)$-orbit type set
$(Q^s)_{(H)}$ can be seen as a manifold.

Consider a topological group $G$ acting continuously on a
topological space $X$ by an action $\Phi : G\times
X\longrightarrow X$. Let $H$ be a subgroup of $G$. As before we
can define the sets $X^H, X_H$, and $X_{(H)}$. Let $N(H):=\{g\in
G\,|\,gHg^{-1}=H\}$ be the normalizer of $H$ in $G$. We can
consider the following well-defined twisted action of $N(H)$ on
$G\times X_H$:
\[
N(H)\times (G\times X_H)\longrightarrow (G\times
X_H),\;\;(h,(g,x))\longmapsto (gh,\Phi_{h^{-1}}(x)).
\]
The orbit space of this free action is denoted by
$G\times_{N(H)}X_H$.

Remark that the continuous map $G\times X_H\longrightarrow
X_{(H)},\;\;(g,x)\longmapsto \Phi_g(x)$ induces a well defined continuous
bijection
\begin{equation}\label{first bijection}
G\times_{N(H)}X_H\longrightarrow
X_{(H)},\;\;[(g,x)]_{N(H)}\longmapsto \Phi_{g}(x),
\end{equation}
where $[(g,x)]_{N(H)}$ denotes the equivalence class of $(g,x)$
relative to the twisted action.

We now consider the continuous map $X_H\longrightarrow
X_{(H)}/G,\;\;x\longmapsto [x]_{G}$, where $[x]_{G}$ denotes the
equivalence class of $x$ in $X_{(H)}$ relative to the action $\Phi
: G\times X_{(H)}\longrightarrow X_{(H)}$. One can show that this
map induces a continuous bijection
\begin{equation}\label{second bijection}
X_H/N(H)\longrightarrow X_{(H)}/G,\;\; [x]_{N(H)}\longmapsto
[x]_{G},
\end{equation}
where $[x]_{N(H)}$ denotes the equivalence class of $x$ in $X_{H}$
relative to the action $\Phi : N(H)\times X_{H}\longrightarrow
X_{H}$.

Remark that $X_H/(N(H)/H)=X_H/N(H)$, and that the action of
$N(H)/H$ on $X_H$ is free.

Applying the result \eqref{first bijection} to the non-smooth
action $\Phi : G\times Q^s\longrightarrow Q^s$ verifying the
hypotheses of Section 2, we obtain a continuous bijection
\[
G\times_{N(H)}Q^s_H\longrightarrow (Q^s)_{(H)}.
\]
Since the twisted action of $N(H)$ on $G\times Q^s_H$ is free and
verifies the hypotheses of Section 2, the orbit space
$G\times_{N(H)}Q^s_H$ is a topological manifold (by Theorem
\ref{TM}) and all the results obtained in the present section are
valid: we can define for $G\times_{N(H)}Q^s_H$ the weak tangent
bundle and the weak differentiable curves. Using the previous
continuous bijection, we can transport all these properties on
$(Q^s)_{(H)}$.

In the same way, applying the result \eqref{second bijection} to
the non-smooth action $\Phi : G\times Q^s\longrightarrow Q^s$
verifying the hypotheses of Section 2, we obtain a continuous
bijection
\[
Q^s_H/(N(H)/H)\longrightarrow (Q^s)_{(H)}/G.
\]
Since the action of $N(H)/H$ on $Q^s_H$ is free and verifies the
hypotheses of Section 2, the orbit space $Q^s_H/(N(H)/H)$ is a
topological manifold and we can define the weak tangent bundle and
the weak differentiable curves on it. Using the previous
continuous bijection, we can transport all these properties on
$(Q^s)_{(H)}/G$.

\section{Differentiable functions on $Q^s/G$}

Now we define a notion of differentiable functions on $Q^s$ and
$Q^s/G$ that is compatible with the notion of
weakly-differentiable curves in $Q^s$ and in $Q^s/G$.

\begin{definition}
\label{weakly differentiable functions}
Let $s> s_0+1$.
\begin{enumerate}
\item[{\rm (i)}] We define the set $C^1_W(Q^s)$ of
\textbf{weakly-differentiable function} on $Q^s$ as follows:
\begin{quote}$f\in C^1_W(Q^s)$ if and only if $f$ is a real-valued
function on $Q^s$ which is $C^1$ relative to the $s-1$
differentiable structure on $Q^s$.\end{quote}

\item[{\rm (ii)}] The set $C^1_W(Q^s/G)$ of
\textbf{weakly-differentiable functions} on $Q^s/G$ is defined as
follows:
\begin{quote}$\varphi\in C^1_W(Q^s/G)$ if and only if $\varphi$
is a real-valued function on $Q^s/G$ such that $\varphi\circ\pi\in
C^1_W(Q^s)$.\end{quote}
\end{enumerate}
\end{definition}

Note that for all $f\in C^1(Q^{s-1})$ we have $f|_{Q^s}\in
C^1_W(Q^s)$. In this sense we can write the inclusion
$C^1(Q^{s-1})\subset C^1_W(Q^s)$. Moreover, for all $f\in
C^1_W(Q^s)$ we can define the function $f\circ j_{(s,s-1)}\in
C^1(Q^s)$. In this sense we can write the inclusion
$C^1_W(Q^s)\subset C^1(Q^s)$. Thus we obtain
\[
C^1(Q^{s-1})\subset C^1_W(Q^s)\subset C^1(Q^s).
\]
If $f\in C^1_W(Q^s)$ and $q\in Q^s$, by definition the
differential
\[
df(q) : T_qQ^s\longrightarrow \mathbb{R}
\]
is a continuous linear map relative to the $s-1$ topology on
$T_qQ^s$. Since $Q^s$ is dense in $Q^{s-1}$, using exponential
charts we obtain that for all $q\in Q^s$, $T_qQ^s$ is dense in
$T_qQ^{s-1}$. Thus we can extend $df(q)$ to the map
\[
d^Wf(q) : T^W_qQ^s\longrightarrow \mathbb{R}
\]
which is linear and continuous on $T^W_qQ^s$ with the topology of
$T_qQ^{s-1}$.

The next result shows that the notion of weakly-differentiable
function in Definition \ref{weakly differentiable functions} is
compatible with the notion of weakly-differentiable curves in Definition
\ref{weakly differentialble curves}. Note that for the first time, we
use the density of the inclusion $j_{(s,s-1)} : Q^s\hookrightarrow
Q^{s-1}$.

\begin{theorem}\label{differential_of_f} Let $s> s_0+1$ and $f\in
C^1_W(Q^s)$. Then for all $q\in Q^s$ and $v_q\in T_q^WQ^s$ we have
\[
d^Wf(q)(v_q)=\left.\frac{d}{dt}\right|_{t=0}f(d(t))
\]
where $d$ is any curve in $C^1_W(I,Q^s)$ such that $d(0)=q$ and
$\dot{d}(0)=v_q$.
\end{theorem}
\textbf{Proof.} Let $d\in C^1_W(I,Q^s)$ be such that
$\dot{d}(0)=v_q$. Since $f\circ d\in C^1(I,\mathbb{R})$ we have
\[
\lim_{t\rightarrow
0}\frac{f(d(t))-f(d(0))}{t}=\left.\frac{d}{dt}\right|_{t=0}f(d(t)).
\]
We shall show that
\[
\lim_{t\rightarrow 0}\frac{f(d(t))-f(d(0))}{t}=d^Wf(q)(v_q).
\]
Using a chart $\varphi^s$ of $Q^s$ at $q$, which is the
restriction of a chart $\varphi^{s-1}$ of $Q^{s-1}$, we have
\begin{align*}
&\frac{f(d(t))-f(d(0))}{t}-d^Wf(q)(v_q)\\
&=\frac{(f\circ(\varphi^s)^{-1})(\varphi^s(d(t)))
-(f\circ(\varphi^s)^{-1})(\varphi^s(d(0)))}{t}\\
&\qquad\qquad-d^W(f\circ(\varphi^s)^{-1})(\varphi^s(d(0)))
\left(T_{d(0)}\varphi^s(\dot{d}(0))\right)\\
&=\frac{(f\circ(\varphi^s)^{-1})(\varphi^s(d(t)))
-(f\circ(\varphi^s)^{-1})(\varphi^s(d(0)))}{t}\\
&\qquad\quad-d(f\circ(\varphi^s)^{-1})(\varphi^s(d(0)))
\left(\frac{\varphi^s(d(t))-\varphi^s(d(0))}{t}
\right)\\
&\qquad\quad+d(f\circ(\varphi^s)^{-1})(\varphi^s(d(0)))
\left(\frac{\varphi^s(d(t))-\varphi^s(d(0))}{t}
\right)\\
&\qquad\quad-d^W(f\circ(\varphi^s)^{-1})(\varphi^s(d(0)))
\left(T_{d(0)}\varphi^s(\dot{d}(0))\right)\\
&=\frac{1}{t}\Big{(}(f\circ(\varphi^s)^{-1})(\varphi^s(d (t)))
-(f\circ(\varphi^s)^{-1})(\varphi^s(d(0)))\\
&\qquad\quad-d(f\circ(\varphi^s)^{-1})(\varphi^s(d(0)))
\left(\varphi^s(d(t))-\varphi^s(d(0))\right)
\Big{)}\\
&\qquad\quad+d^W(f\circ(\varphi^s)^{-1})(\varphi^s(d(0)))
\left(\frac{\varphi^s(d(t))-\varphi^s(d(0))}
{t}-T_{d(0)}\varphi^s(\dot{d}(0))\right).
\end{align*}

By continuity of the linear map
$d^W(f\circ(\varphi^s)^{-1})(\varphi^s(d(0)))$ with respect to the
$s-1$ norm, the last term converges to $0$. By
weak-differentiability of $f\circ(\varphi^s)^{-1}$, for all
$\varepsilon >0$ we can choose $\delta>0$ such that if
$|t|<\delta$ the first term is less than
\[
\varepsilon\left\|\frac{\varphi^s(d(t))-\varphi^s(d(0))}{t}\right\|_{s-1}.
\]
Since $\varphi^s\circ d\in C^1(I,T_qQ^{s-1})$, this expression
 converges to
\[
\varepsilon\left\|\left.\frac{d}{dt}\right|_{t=0}\varphi^s(d
(t))\right\|_{s-1}.
\]
Thus we obtain that $\frac{f(d(t))-f(d(0))}{t}$ converges to
$d^Wf(q)(v_q)$.\;\;\;\;$\blacksquare$

\medskip

Inspired by the previous result we now define the differential of
a function in $C^1_W(Q^s/G)$.

\begin{definition} Let $s> s_0+2$ and $\varphi\in
C^1_W(Q^s/G)$. For $\omega\in Q^s/G$ and $v_\omega\in
T^W_\omega(Q^s/G)$, the \textbf{differential} of $\varphi$ at
$\omega$ in direction $v_\omega$ is defined by
\[
d^W\varphi(\omega)(v_\omega):=\left.\frac{d}{dt}\right|_{t=0}\varphi(c(t)),
\]
where $c\in C^1_W(I,Q^s/G)$ is any curve such that $c(0)=\omega$
and $\dot{c}(0)=v_\omega$.
\end{definition}

\begin{theorem}\label{Horizontalchainrule} Let $s> s_0+2$ and $\varphi\in
C^1_W
(Q^s/G)$. Then
the differential of $\varphi$ at $\omega\in Q^s/G$ in direction
$v_\omega\in T^W_\omega(Q^s/G)$ is well-defined, that is, it does
not depend on the choice of the curve $c\in C^1_W(I,Q^s/G)$ satisfying
$c(0)=\omega$ and $\dot{c}(0)=v_\omega$.
Moreover, for $f:=\varphi\circ\pi\in C^1_W(Q^s)$ and
$q_\omega\in\pi^{-1}(\omega)$ we have
\[
d^W\varphi(\omega)(v_\omega)=d^Wf(q_\omega)(\operatorname{Hor}_{q_\omega}
(v_\omega))
\]
\end{theorem}
\textbf{Proof.} Let $c\in C^1_W(I,Q^s/G)$ satisfy
$\dot{c}(0)=v_\omega$  such that
$c=\pi\circ d \in C ^1_W(I', Q ^s)$, $I' \subset I$ (Definition
\ref{weakly differentialble curves} of
$C^1_W(I,Q^s/G)$). Note that we have $\varphi\circ c=\varphi\circ
\pi\circ d\in C^1(I,\mathbb{R})$ since $\varphi\circ\pi\in
C^1_W(Q^s)$. So we can compute the derivative of $\varphi\circ c$.
With $f:=\varphi\circ\pi$ we have
\begin{align*}
\left.\frac{d}{dt}\right|_{t=0}\varphi(c(t))
&=\left.\frac{d}{dt}\right|_{t=0}f (d(t))\\
&=d^Wf(d(0))(\dot{d}(0))\,\text{ by Theorem \ref{differential_of_f}}\\
&=d^Wf(d(0))(\operatorname{hor}_{d(0)}(\dot{d}(0)))\\
&=d^Wf(d(0))(\operatorname{Hor}_{d(0)}(T^W_{d(0)}\pi(\dot{d}(0)))\\
&=d^Wf(d(0))(\operatorname{Hor}_{d(0)}(\dot{c}(0)))\,\text{ by
Definition
\ref{definition_of_weak} (v)}\\
&=d^Wf(q_\omega)(\operatorname{Hor}_{q_\omega}(v_\omega)).
\end{align*}
The third equality follows from the fact that $d^Wf(q_\omega)$ vanishes
on the vertical subspace. Namely, since $f=\varphi\circ\pi$ is
$G$-invariant, that is, $f\circ\Phi_g=f$ for all $g\in G$, we obtain
\begin{equation}
\label{weak derivative horizontal}
0=\left.\frac{d}{dt}\right|_{t=0}f(\Phi_{\operatorname{exp}(t\xi)}(q))
=d^Wf(q)\left(\left.\frac{d}{dt}\right|_{t=0}\Phi_{\operatorname{exp}(t\xi)}(q)
\right)
=d^Wf(q)(\xi_{Q^s}(q)),
\end{equation}
where we used that $t\longmapsto
\Phi_{\operatorname{exp}(t\xi)}(q)$ is in
$C^1_W(I,Q^s)$.\;\;\;\;$\blacksquare$

\section{Application to fluid dynamics}

We consider the motion of an incompressible ideal fluid in a
compact oriented Riemannian manifold $M$ with boundary. It is well
known that the configuration space is $\mathcal{D}_\mu^s (M),
s>\frac{\operatorname{dim}(M)}{2}+1$, the Hilbert manifold of
volume preserving $H^s$-diffeomorphisms of $M$, and that the
appropriate Lagrangian is given by the weak $L^2$
Riemannian metric
\[
\langle\!\langle u_\eta,v_\eta
\rangle\!\rangle_\eta=\int_Mg(\eta(x))(u_\eta(x),v_\eta(x))\mu
(x),\,\,\,u_\eta,v_\eta\in
T_\eta \mathcal{D}^s_\mu(M),
\]
where $g$ is the Riemannian metric on $M$ and $\mu$ is the volume
form induced by $g$. This Lagrangian is
invariant under the two following commuting actions:
\[
R : \mathcal{D}_\mu^s(M)\times T\mathcal{D}_\mu^s(M)
\longrightarrow
T\mathcal{D}_\mu^s(M),\,\,R(\eta,v_\xi)=R_\eta(v_\xi):=v_\xi\circ
\eta
\]
\[
L : Iso^+\times T\mathcal{D}_\mu^s(M)\longrightarrow T\mathcal{D}
_\mu^s(M),\,\,L(i,v_\xi)=L_i(v_\xi):=Ti\circ v_\xi,
\]
where $Iso^+:=Iso^+(M,g)$ denotes the group of Riemannian
isometries of $(M,g)$ which preserve the orientation.
Since $M$ is compact, it follows that $Iso^+$ is a
compact Lie group of dimension $\leq\frac{n(n-1)}{2},\,
n=\operatorname{dim}(M)$.

We denote by $\mathfrak{iso}^+:=T_eIso^+$ the Lie algebra of
$Iso^+$. From Corollary $5.4$ of \cite{Eb1968}, we know that
$Iso^+$ can be seen as a submanifold of
$\mathcal{D}^r(M),\,r>\frac{\operatorname{dim}(M)}{2}+1$. The
tangent space at the identity of $Iso^+$, viewed as a submanifold
of $\mathcal{D}^r(M)$, consists of smooth vector fields on $M$
whose flows are curves in $Iso^+$, that is, it consists of the
Killing vector fields:
\[
\mathfrak{X}_K(M)=\{X\in\mathfrak{X}(M)\,|\,L_Xg=0\},
\]
where $L_Xg$ is the Lie derivative of $g$ along $X$. The
correspondence between $\mathfrak{iso}^+$ and $\mathfrak{X}_K(M)$
is given by:
\[
\mathfrak{iso}^+\longrightarrow
\mathfrak{X}_K(M),\;\;\xi\longmapsto X_\xi,
\]
where $X_\xi(x):=T_eEv_x(\xi)$ and $Ev_x: Iso^+\longrightarrow
M,\;\;Ev_x(i):=i(x)$ is the evaluation map at $x$. Indeed, one can
see that the flow of $X_\xi$ is given by
$\operatorname{exp}(t\xi)$, where $\operatorname{exp}$ is the
exponential map of $Iso^+$.

Note that for all $i\in Iso^+$ we have $i_*\mu=\mu$, where $\mu$
is the volume form associated to $g$, so $i$ is volume preserving
and $X_\xi$ is divergence free. For example, setting
$M:=\mathbb{B}^n$, the closed unit ball in $\mathbb{R}^n$, we have
$Iso^+=$ SO(n) and for all $\xi\in\mathfrak{so(n)}$ we have
\[
X_\xi(x)=T_eEv_x(\xi)=\left.\frac{d}{dt}\right|_{t=0}\operatorname{exp}(t\xi)
(x)=\xi
x.
\]

Our goal is to carry out the Poisson reduction by stages
associated to these two commuting actions. For finite dimensional
manifolds and Lie groups the reduction by stages procedure (see
\cite{MaMiOrPeRa2006}) guarantees that the two stage reduction by
the two commuting group actions yields the same result as the
reduction by the product group. In our case this one step
reduction by the product group $\mathcal{D}_\mu^s(M) \times Iso^+$
cannot be carried out because $\mathcal{D}_\mu^s(M)$ is not a Lie
group and the action is not smooth in the usual sense. However, we
shall see that the two step reduction, first by
$\mathcal{D}_\mu^s(M)$ and then by the compact Lie group $Iso^+$,
can be carried out in view of the general results proved above.

The result we shall obtain is a non-smooth generalization of the
following theorem for proper smooth Lie group actions on finite
dimensional manifolds in the physically relevant case of the Euler
equations. Let $G \times M \rightarrow M$ be a smooth proper
action of the Lie group $G$ on the Poisson manifold $(M,
\{\,,\})$. If $(H)$ is an orbit type then $M_{(H)}/G$ is a smooth
Poisson manifold. The Poisson structure is obtained by push
forward of the natural quotient Poisson bracket on $M_H/N(H)$,
where $N(H)$ is the normalizer of $H$ in $G$. For details see
\cite{FeOrRa2006}. In our case we shall proceed in the following
way. First we reduce $T \mathcal{D}^s_\mu(M)$ by the right action
of $\mathcal{D}^s_\mu(M)$ and obtain $T\mathcal{D}_\mu^s(M)/
\mathcal{D}_\mu^s(M) = \mathfrak{X}^s_{div}(M)$. The action of
$Iso^+$ drops to a Poisson action on $\mathfrak{X}^s_{div}(M)$.
Then we find the explicit Poisson bracket on the isotropy type
manifolds $\mathfrak{X}^s_{div}(M)_H$ and on its quotient
$\mathfrak{X}^s_{div}(M)_H/N(H)$, for any isotropy subgroup $H
\subset Iso^+$. We close by presenting the reduced Euler equations
on this quotient and discuss in what sense the flow is Poisson.

\subsection{Reduction by
$\mathcal{D}^s_\mu(M)$}

Reduction by $\mathcal{D}_\mu^s(M)$ is well known (see
\cite{EbMa1970}) and leads to the Euler equations for an
ideal incompressible fluid on the first reduced space
$\mathfrak{X}^s_{div}(M)= T\mathcal{D}_\mu^s(M)/
\mathcal{D}_\mu^s(M)$ consisting of $H^s$ divergence free vector
fields on $M$ that are tangent to the boundary. The most
fundamental fact is the existence of the smooth geodesic spray
$\mathcal{S}\in\mathfrak{X}^{C^\infty}(T\mathcal{D}^s_\mu(M))$ of
the weak Riemannian manifold
$(\mathcal{D}^s_\mu(M),\langle\!\langle\,, \rangle\!\rangle)$. The
following reduction theorem can by found in \cite{EbMa1970}.

\begin{theorem} Let $\eta(t)\subset
\mathcal{D}^s_\mu(M),s>\frac{\operatorname
{dim}(M)}{2}+1$,
be a curve in $\mathcal{D}^s_\mu(M)$ and let
$u(t):=R_{\eta(t)^{-1}}(\dot{\eta}(t))=\dot{\eta}(t)\circ\eta(t)^{-1}
\in\mathfrak{X}^s_{div}(M)$. Then the following properties are
equivalent.
\begin{enumerate}
\item[{\rm (i)}] $\eta(t)$ is a geodesic of
$(\mathcal{D}^s_\mu(M),\langle\!\langle\,, \rangle\!\rangle)$.
\item[{\rm (ii)}] $V(t):=\dot\eta(t)$ is a solution of $\dot
V(t)=\mathcal{S}(V(t))$.
 \item[{\rm (iii)}] $u(t)$
is a solution of the Euler equations
\[
\partial_tu(t)+\nabla_{u(t)}u(t)=-\operatorname{grad}p(t)
\]
for some scalar function $p(t) : M\longrightarrow \mathbb{R}$
called the pressure. \end{enumerate} Moreover the solution $u$ of
the Euler equation is in $C^0(I,\mathfrak{X}^s_{div}(M))\cap
C^1(I,\mathfrak{X}^{s-1}_{div}(M))$
\end{theorem}

Note the the Euler equations can be written as
\[
\partial_tu(t)+P_e(\nabla_{u(t)}u(t))=0
\]
where $P_e$ denotes the projection on the first factor of the
Hodge decomposition
\[
\mathfrak{X}^r(M)=\mathfrak{X}^r_{div}(M)\oplus\operatorname
{grad}(H^{r+1}(M)),\quad r\geq 0.
\]
Denoting by $\pi_R : T\mathcal{D}^s_\mu(M)\longrightarrow
\mathcal{D}^s_\mu(M),\;\;\pi_R(u_\eta):=u_\eta\circ\eta^{-1}$, the
projection associated to the reduction by $\mathcal{D}^s_\mu(M)$
we obtain the following commutative diagram
$$\xymatrix{
T\mathcal{D}^s_\mu(M) \ar[r]^{F_t} \ar[d]_{\pi_R}&
T\mathcal{D}^s_\mu(M) \ar [d]^{\pi_R}\\
\mathfrak{X}^s_{div}(M)\ar[r]^{\widetilde{F}_t} &
\mathfrak{X}^s_{div}(M).\\
}$$ where $F_t$ is the flow of $\mathcal{S}$ and $\widetilde{F}_t$
is the flow of the Euler equations. Formally,
all these maps are Poisson, as is the case in the standard
Poisson reduction procedure. However, since our manifolds are
infinite dimensional, some difficulties arise. First, the
symplectic form on $T\mathcal{D}^s_\mu(M)$ is only weak since the
Lagrangian is given by a $L^2$ metric. Second,
$\mathcal{D}^s_\mu(M)$ is not a Lie group since left
multiplication and inversion are not smooth. \cite{VaMa2005} have
resolved these difficulties by carefully analyzing the function
spaces on which Poisson brackets are defined and carrying out a
non-smooth Lie-Poisson reduction that takes into account all
analytical difficulties. We recall below some results of this
paper that we will use later on.

\medskip

\noindent\textbf{Poisson brackets on $T\mathcal{D}^s_\mu(M)$}
\medskip

For $F : T\mathcal{D}^s_\mu(M)\longrightarrow \mathbb{R}$ of class
$C^1$, we define the horizontal partial derivative of $F$ by
\[
\frac{\partial F}{\partial \eta} : T \mathcal{D}^s_\mu(M)
\rightarrow T ^\ast \mathcal{D}^s_\mu(M)
\]
such that
\[
\frac{\partial F}{\partial \eta}(u_\eta)(v_\eta) : =
\left.\frac{d}{dt}\right|_{t=0}F(\gamma(t)),
\]
where $\gamma(t) \subset T \mathcal{D}^s_\mu(M)$ is a smooth path
defined in a neighborhood of zero, with base point denoted by
$\eta(t) \subset \mathcal{D}^s_\mu(M)$, satisfying the following
conditions:
\begin{itemize}
\item $\gamma(0) = u_\eta$ \item $\dot{\eta}(0) = v_\eta$ \item
$\gamma$ is parallel, that is, its  covariant derivative
associated to the metric $\langle\!\langle\,, \rangle\!\rangle$
vanishes.
\end{itemize}
The vertical partial derivative
\[
\frac{\partial F}{\partial u} : T \mathcal{D}^s_\mu(M) \rightarrow
T ^\ast \mathcal{D}^s_\mu(M)
\]
of $F $ is defined as the usual fiber derivative, that is,
\[
\frac{\partial F}{\partial u}(u_\eta)(v_\eta) : =
\left.\frac{d}{dt}\right|_{t=0} F(u_\eta + t v_\eta).
\]
These derivatives naturally induce corresponding functional
derivatives relative to the weak Riemannian metric
$\langle\!\langle\,, \rangle\!\rangle$. The horizontal and
vertical functional derivatives
\[
\frac{\delta F}{\delta \eta }, \frac{\delta F}{\delta u} : T
\mathcal{D}^s_\mu(M) \rightarrow T \mathcal{D}^s_\mu(M)
\]
are defined by the equalities
\[
\left\langle\!\!\!\left\langle\frac{\delta F}{\delta
\eta}(u_\eta),v_\eta \right\rangle\!\!\!\right\rangle =
\frac{\partial F}{\partial \eta}(u_\eta)(v_\eta) \quad \text{ and
} \quad \left\langle\!\!\!\left\langle\frac{\delta F}{\delta
u}(u_\eta),v_\eta \right\rangle\!\!\!\right\rangle =
\frac{\partial F}{\partial u}(u_\eta)(v_\eta)
\]
for any $u_\eta, v_ \eta\in T\mathcal{D}^s_\mu(M)$. Note that due
to the weak character of $\langle\!\langle\,, \rangle\!\rangle$,
the existence of the functional derivatives is not guaranteed. But
if they exist, they are unique.

We define, for $k\geq 1$ and
$r,t>\frac{\operatorname{dim}(M)}{2}+1$ :
\[
C^k_r(T\mathcal{D}^t_\mu(M)):=\left\{F\in
C^k(T\mathcal{D}^t_\mu(M)) \Big|\exists\,\frac{\delta F}{\delta
\eta},\frac{\delta F}{\delta u} : T\mathcal{D}
^t_\mu(M)\longrightarrow T\mathcal{D}^r_\mu(M)\right\}.
\]

With these definitions the Poisson bracket of $F,G\in
C^k_r(T\mathcal{D}^t_\mu(M))$ is given by
\begin{equation}
\label{poisson bracket} \{F,G\}(u_\eta)
=\left\langle\!\!\!\left\langle \frac{\delta F}{\delta \eta}
(u_\eta),\frac{\delta G}{\delta
u}(u_\eta)\right\rangle\!\!\!\right\rangle -
\left\langle\!\!\!\left\langle\frac{\delta F}{\delta u}(u_\eta),
\frac{\delta G}{\delta
\eta}(u_\eta)\right\rangle\!\!\!\right\rangle.
\end{equation}

\noindent\textbf{Poisson brackets on $\mathfrak{X}^s_{div}(M)$}
\medskip

For $k\geq 1$ and $r\geq 0$ we define the set
\[
C^k_{r}(\mathfrak{X}^s_{div}(M)) :=\left\{f\in
C^k(\mathfrak{X}^s_{div}(M)) |\exists\,\delta f
:\mathfrak{X}^s_{div}(M)\longrightarrow
\mathfrak{X}^r_{div}(M)\right\}
\]
where $\delta f$ is the functional derivative of $f$ with respect
to the inner product $\langle \,\,,\,\rangle:=\langle\!\langle\,,
\rangle\!\rangle_{id}$, that is
\[
\langle\delta f(u),v\rangle=Df(u)(v),\quad
\forall\,u,v\in\mathfrak{X}^s_{div}(M).
\]
For $k\geq 1$ and $r>\frac{\operatorname{dim}(M)}{2}+1$, the
Poisson bracket of $f,g\in C^k_r(\mathfrak{X}^s_{div}(M))$ is
defined by
\begin{equation}
\label{Lie-Poisson bracket} \{f,g\}_+(u):=\langle u,[\delta
g(u),\delta f(u)]\rangle,\quad \forall
u\in\mathfrak{X}^s_{div}(M).
\end{equation}

The next theorem summarizes the principal results of
\cite{VaMa2005}.

\begin{theorem} \label{VaMa} Let $k\geq 1$.
\begin{enumerate}\item[{\rm (i)}]Let $F_t$ be the flow of the
geodesic spray $\mathcal{S}$ and let $t_1\geq
t_2>\frac{\operatorname{dim}(M)}{2}+1$. Then for all $G,H\in
C^k_{t_2}(T\mathcal{D}^{t_1}_\mu(M))$ we have
\[
\{G\circ F_t,H\circ F_t\}=\{G,F\}\circ F_t
\]
on $T\mathcal{D}^{t_1}_ {\mu,D}$.

\item[{\rm (ii)}]Let $r>\frac{\operatorname{dim}(M)}{2}+1$ satisfy
$s+k\geq r$. Then for all $f,g\in
C^k_r(\mathfrak{X}^s_{div}(M))$ we have
\[
\{f\circ\pi_R,g\circ\pi_R\}(u_\eta)
=\left(\{f,g\}_+\circ\pi_R\right) (u_\eta),\quad \forall u_\eta\in
T\mathcal{D}^{s+k}_\mu(M).
\]
\item[{\rm (iii)}] Let $\widetilde{F}_t$ be the flow of the Euler
equations and let $r>2$ be such that $s+k\geq r$. Then for all $f,g\in
C^k_r(\mathfrak{X}^s_{div}(M))$ we have
\[
\{f\circ\tilde{F}_t,g\circ\tilde{F}_t\}_+(u)
=\left(\{f,g\}_+\circ\tilde{F}_t\right) (u),\quad \forall
u\in\mathfrak{X}^{s+2k}_{div}(M).
\]
\end{enumerate}
\end{theorem}

In the next subsections we will carry out in a precise sense the
second stage reduction, that is, the reduction by the group
$Iso^+$.

\subsection{Action of $Iso^+$ on $\mathfrak{X}^s_{div}(M)$}

Recall that the action of $Iso^+$ on the tangent bundle of
$\mathcal{D}^s_\mu(M)$ is
\[
L : Iso^+\times T\mathcal{D}_\mu^s(M)\longrightarrow T\mathcal{D}
_\mu^s(M),\,\,L_i(v_\xi)=Ti\circ v_\xi.
\]
Since $R$ and $L$ commute, $L$ induces the action
\[
l : Iso^+\times\mathfrak{X}^s_{div}(M) \longrightarrow
\mathfrak{X}^s_{div}(M),\;\;l_i(u)=i_*u.
\]
Indeed, $l_i(u):=\pi_R(L_i(u))=L_i(u)\circ i^{-1}=Ti\circ u\circ
i^{-1}=i_*u$.

Remark that the action $l$ is not smooth on
$Iso^+\times\mathfrak{X}^s_{div}(M)$. However, as we shall see
below, $l$ verifies all the hypothesis in Section 2. Thus
the general theory developed in Sections $3$ to $5$ is directly
applicable to the present case.

In general, the action $l$ is not free. Indeed, consider the
particular case $M=\mathbb{D}\subset\mathbb{R}^2$, the closed unit
disc, so we have $Iso^+=$ SO(2) and one sees that the vector field
\[
z : \mathbb{D}\longrightarrow \mathbb{R}^2, \qquad  z(a,b)=(-b,a)
\]
is divergence free, tangent to the boundary, with isotropy group
equal to SO(2).

We now proceed to the verification of the hypotheses of the
Section 2. Consider the collection
$\left\{\mathfrak{X}^s_{div}(M)\,|\,s>\frac{\operatorname{dim}(M)}{2}\right\}$
of Hilbert spaces. It is clear that the inclusions $j_{(r,s)}$ are
smooth with dense range as well as their tangent maps. We endow
$\mathfrak{X}^s_{div}(M)$ with the $L^2$ inner product
\[
\langle u,v\rangle:=\int_Mg(x)(u(x),v(x))\mu(x)
\]
which is the value at the identity of the weak Riemannian
metric $\langle\!\langle\,,\rangle\!\rangle$. Note that for all
$i\in Iso^+$, we have
\[
\langle l_i(u),l_i(v)\rangle=\langle u,v\rangle.
\]

By \cite{Eb1968} we know that the map
\[
\mathcal{D}^r(M)\times\mathfrak{X}^s_{div}(M)\longrightarrow
\mathfrak{X}^s_{div}(M),\;\;(i,u)\longmapsto i_*u
\]
is continuous for $r$ sufficiently large;  in \cite{Eb1968},
$\mathcal{D}^r(M)$ acts by pull back on Riemannian metrics and
here on vector fields, but the proofs are similar. Since
$Iso^+$ can be seen as a smooth submanifold of
$\mathcal{D}^r(M),\,r>\frac{\operatorname{dim}(M)}{2}+1$, we
obtain that
\[
Iso^+\times\mathfrak{X}^s_{div}(M)\longrightarrow
\mathfrak{X}^s_{div}(M),\;\;(i,u)\longmapsto i_*u
\]
is continuous. Since $Iso^+$ is a compact Lie group, the action is
proper.

We now check that $l$ verifies the hypotheses (\ref{phi_g_smooth})
and (\ref{phi_c_one}). It is obvious that for all $i\in Iso^+$ the
map
\[
l_i : \mathfrak{X}^s_{div}(M)\longrightarrow
\mathfrak{X}^s_{div}(M)\;\;l_i(u)=i_*u
\]
is smooth since it is linear and continuous. It remains to show
that for all $u\in\mathfrak{X}^s_{div}(M),
s>\frac{\operatorname{dim}(M)}{2}+1$ the map
\[
l^u : Iso^+\longrightarrow \mathfrak{X}^{s-1}_{div}(M),\,\,
l^u(i):=i_*u
\]
is $C^1$. By the proof of Proposition 3.4 of \cite{Eb1968}, the
map
\[
\mathcal{D}^r(M)\longrightarrow H^{s-1}(M,TM),\;\;\eta\longmapsto
u\circ\eta
\]
is $C^1$ for $r$ sufficiently large. Since $Iso^+$ is a smooth
submanifold of $\mathcal{D}^r(M)$, the map
\[
Iso^+\longrightarrow H^{s-1}(M,TM),\;\;i\longmapsto u\circ i
\]
is $C^1$. Using that
\[
Iso^+\longrightarrow H^s(M,TM),\,\,i\longmapsto Ti\circ u
\]
and $i\longmapsto i^{-1}$ are smooth, we obtain the desired
result.
\medskip

The infinitesimal generator associated to $\xi\in
\mathfrak{iso}^+$ is given by
\begin{equation}\label{infinitesimal generators}
\xi_{\mathfrak{X}^s_{div}(M)}(u)=
\left.\frac{d}{dt}\right|_{t=0}l_{\operatorname{exp}(t\xi)}(u)
=\left.\frac{d}{dt}\right|_{t=0}\operatorname{exp}(t\xi)_*u
=-L_{X_\xi}u=[u,X_\xi],
\end{equation}
where $X_\xi(x)=T_eEv_x(\xi)$ is the Killing vector field
generated by the flow $\operatorname{exp}(t\xi)$ and $[\,,]$ is
the Jacobi-Lie bracket of vector fields. Remark that we have, as
expected, $[u,X_\xi]\in\mathfrak{X}^{s-1}_{div}(M)$, since the
Jacobi-Lie bracket of divergence free vector fields remains
divergence free.

In the particular case $M=\mathbb{B}^n$ and $Iso^+=$ SO(n) we have
$l_A(u)(x)=A u(A^{-1}x)$, so we compute directly
\[
\xi_{\mathfrak{X}^s_{div}(\mathbb{B}^n)}(u)(x)
=\left.\frac{d}{dt}\right|_{t=0}
\operatorname{exp}(t\xi)u(\operatorname{exp}(-t\xi)x)=\xi
u(x)-Du(x)(\xi x).
\]
Since $X_\xi(x)=\xi x$ and $DX_\xi(x)=\xi$ we obtain
\[
\xi u(x)-Du(x)(\xi
x)=DX_\xi(x)(u(x))-Du(x)(X_\xi(x))=[u,X_\xi](x),
\]
so we recover the expression \eqref{infinitesimal generators}.

\medskip

Now we can use the results in Section $3$ to $5$ with
$Q^s=\mathfrak{X}^s_{div}(M)$ and
$s_0=\frac{\operatorname{dim}(M)}{2}$. We obtain the following
results:

\begin{enumerate}

\item[{\rm (i)}] \textit{For $s>\frac{\operatorname{dim}(M)}{2}+1$
and $H\subset Iso^+$ a closed subgroup,
$\mathfrak{X}^s_{div}(M)_H$ is an open set in the vector subspace
$\mathfrak{X}^s_{div}(M)^H$ of $\mathfrak{X}^s_{div}(M)$. As we
did in the general case, we can define the weak tangent bundle
$T^W\mathfrak{X}^s_{div}(M)$. Since $\mathfrak{X}^s_{div}(M)$ is a
vector space, we use the notation
$T^W\mathfrak{X}^s_{div}(M)=\mathfrak{X}^s_{div}(M)\times
\mathfrak{X}^s_{div}(M)^W$.}

\item[{\rm (ii)}] \textit{Since the action $N(H)/H\times
\mathfrak{X}^s_{div}(M)_H\longrightarrow
\mathfrak{X}^s_{div}(M)_H$ is free and verifies the hypotheses in
Section 2, we obtain that for
$s>\frac{\operatorname{dim}(M)}{2}+1$,
$\mathfrak{X}^s_{div}(M)_H/N(H)=\mathfrak{X}^s_{div}(M)_H/(N(H)/H)$
is a topological manifold.}

\item[{\rm (iii)}] \textit{For
$s>\frac{\operatorname{dim}(M)}{2}+2$ we can define the weak
tangent bundles
\[
T^W(\mathfrak{X}^s_{div}(M)),\; T^W(\mathfrak{X}^s_{div}(M)_H), \;
\text{ and }\;T^W(\mathfrak{X}^s_{div}(M)_H/N(H)),
\]
the spaces of weakly-differentiable curves
\[
C^1_W(I,\mathfrak{X}^s_{div}(M)),\;
C^1_W(I,\mathfrak{X}^s_{div}(M)_H),\;\text{ and }\;
C^1_W(I,\mathfrak{X}^s_{div}(M)_H/N(H)),
\]
the spaces of weakly-differentiable functions
\[
C^1_W(\mathfrak{X}^s_{div}(M)),\;
C^1_W(\mathfrak{X}^s_{div}(M)_H),\;\text{ and }\;
C^1_W(\mathfrak{X}^s_{div}(M)_H/N(H)),
\]
the weak-tangent map
$T^W\pi_H:T^W\mathfrak{X}^s_{div}(M)_H\longrightarrow
T^W(\mathfrak{X}^s_{div}(M)_H/N(H))$, where $\pi_H :
\mathfrak{X}^s_{div}(M)_H\longrightarrow
\mathfrak{X}^s_{div}(M)_H/N(H)$ is the orbit map, the
horizontal lift
\[
\operatorname{Hor}_{u_\omega}^H :
T^W_\omega(\mathfrak{X}^s_{div}(M)_H/N(H))\longrightarrow
T^W_{u_\omega}(\mathfrak{X}^s_{div}(M)_H),
\]
and the vertical and horizontal projections $\operatorname{hor}_u^H$
and $\operatorname{ver}_u^H$ associated to the decomposition}
\[
T^W_u(\mathfrak{X}^s_{div}(M)_H)=H^W_u(\mathfrak{X}^s_{div}(M)_H)\oplus
V^W_u(\mathfrak{X}^s_{div}(M)_H).
\]

\end{enumerate}

\subsection{Dynamics on $(\mathfrak{X}^s_{div}(M),\{\,,\}_+)$}

In this subsection we will show in which sense the action of
$Iso^+$ is canonical with respect to the Lie-Poisson bracket
$\{\,,\}_+$. We will prove existence and uniqueness of the
Hamiltonian vector field associated to a Hamiltonian $h\in
C^1_r(\mathfrak{X}^s_{div}(M))$ and then we will show that the
Hamiltonian vector field associated to the Euler equations takes
values in $T^W(\mathfrak{X}^s_{div}(M))$. Finally, we will see
that the law of conservation of the isotropy is still valid in our
case.

\begin{lemma} \label{Iso^+invariance} Let $i\in Iso^+$.
\begin{enumerate}
\item[{\rm (i)}] For all $f\in
C^1_r(\mathfrak{X}^s_{div}(M)),r,s\geq 0$ we have $f\circ l_i\in
C^1_r(\mathfrak{X}^s_{div}(M))$ and
\[
\delta (f\circ l_i)=l_i^{-1}\circ\delta f\circ l_i.
\]
\item[{\rm (ii)}] For all $u,v\in\mathfrak{X}^s_{div}(M),
s>\frac{\operatorname{dim}(M)}{2}+1$, we have
\[
[l_i(u),l_i(v)]=l_i([u,v]).
\]
\item[{\rm (iii)}] The Hodge projector $P_e :
\mathfrak{X}^s(M)\longrightarrow \mathfrak{X}^s_{div}(M), s\geq 0$,
is $Iso^+$-invariant, that is,
\[
P_e\circ l_i=l_i\circ P_e.
\]
\end{enumerate}
\end{lemma}
\textbf{Proof.} (i) By the chain rule and using that $f\in
C^1_r(\mathfrak{X}^s_{div}(M))$ we have
\begin{align*} D(f\circ l_i)(u)(v)
&=Df(l_i(u))(Dl_i(u)(v))\\
&=Df(l_i(u))(l_i(v))\\
&=\langle \delta f(l_i(u)),l_i(v)\rangle\\
&=\langle l_i^{-1}(\delta f(l_i(u))),v\rangle.
\end{align*}\\
(ii) This is a consequence of the relation $i_*[u,v]=[i_*u,i_*v]$.\\
(iii) Decompose $u\in\mathfrak{X}^s_{div}(M)$ as
$u=P_e(u)+\operatorname{grad} f$, so for all $i\in Iso^+$ we
obtain $l_i(u)=l_i(P_e(u))+l_i(\operatorname{grad} f)$. A direct
computation using the chain rule gives the equality
$l_i(\operatorname{grad} f)=\operatorname{grad}(f\circ i^{-1})$.
Thus we can write $l_i(u)=l_i(P_e(u))+\operatorname{grad}(f\circ
i^{-1})$. Using that the Hodge decomposition is unique gives
$P_e(l_i(u))=l_i(P_e(u))$.\;\;\;\;$\blacksquare$

\medskip

The next theorem shows that the $Iso^+$-action is Poisson relative
to $\{\,,\}_+$.

\begin{theorem} For all $i\in Iso^+$ and $f,g\in C^1_r(\mathfrak{X}^s_{div}
(M)), r,s>\frac{\operatorname{dim}(M)}{2}+1$, we have
\[
\{f\circ l_i,g\circ l_i\}_+=\{f,g\}_+\circ l_i
\]
\end{theorem}
\textbf{Proof.} By Lemma \ref{Iso^+invariance} and
$Iso^+$-invariance of $\langle\,,\rangle$ we have
\begin{align*}\{f\circ l_i,g\circ l_i\}_+(u)&=\langle u,[\delta
(g\circ l_i)(u),\delta (f\circ l_i)(u)]\rangle\\
&=\langle u,[(l_i^{-1}\circ\delta g\circ
l_i)(u),(l_i^{-1}\circ\delta g\circ
l_i)(u)]\rangle\\
&=\langle u,l_i^{-1}([\delta g(l_i(u)),\delta f(l_i(u))])\rangle\\
&=\langle l_i(u),[\delta g(l_i(u)),\delta f(l_i(u))]\rangle\\
&=(\{f,g\}_+\circ l_i)(u).\;\;\;\;\blacksquare
\end{align*}

\begin{theorem}\label{Poissonformulation} Let
$h\in C^1_r(\mathfrak{X}^s_{div}(M)),
r,s>\frac{\operatorname{dim}(M)}{2}+1$. Then there exists a unique
Hamiltonian vector field $X_h$ such that
\[
Df(u)(X_h(u))=\{f,h\}_+(u),\,\,\text{ for all } f\in
C^1_r(\mathfrak{X}^s_{div}(M)).
\]
Moreover $X_h$ is given by
\[
X_h(u)=-P_e(\nabla_{\delta h(u)}u+\nabla\delta h(u)^T\!\!\cdot u),
\]
where $\nabla$ denotes the Levi-Civita covariant derivative and
the upper index $T$ denotes the transpose with respect to the
Riemannian metric $g$. In particular, for $h(u)=\frac{1}{2}\langle
u,u\rangle$ we have
\[
X_h(u)=-P_e(\nabla_uu)
\]
and $X_h : \mathfrak{X}^s_{div}(M)\longrightarrow
T^W\mathfrak{X}^s_{div}(M)$. As usual, we think of the vector $X_h
(u)$ also as $X_h(u)=(u,-P_e(\nabla_uu))$ when the need arises.

\end{theorem}
\textbf{Proof.} For $f,h\in C^1_r(\mathfrak{X}^s_{div}(M))$
integration by parts in the first term gives
\begin{align*}
\{f,h\}_+(u)&=\langle u,[\delta h(u),\delta f(u)]\rangle\\
&=\langle u,\nabla_{\delta h(u)}\delta f(u) - \nabla_{\delta f(u)}\delta
h(u)
\rangle\\
&=-\langle \nabla_{\delta h(u)}u,\delta f(u)\rangle-\langle
\nabla\delta h(u)^T\!\!\cdot u,\delta f(u)\rangle\\
&=-\langle \delta f(u),\nabla_{\delta h(u)}u+\nabla\delta h(u)^T\!\!\cdot
u\rangle\\
&=-\langle \delta f(u),P_e(\nabla_{\delta h(u)}u+\nabla\delta
h(u)^T\!\!\cdot u)\rangle.
\end{align*}
Using that $Df(u)(X_h(u))=\langle \delta f(u),X_h(u)\rangle$ we
obtain by density the existence and uniqueness of
$X_h(u)=-P_e(\nabla_{\delta h(u)}u+\nabla\delta h(u)^T\!\!\cdot
u)$.

In particular, for the reduced Hamiltonian $h(u)=\frac{1}{2}\langle
u,u\rangle$ we have $h\in C^1_s(\mathfrak{X}^s_{div}(M))$ and
$\delta h(u)=u$. Thus $X_h(u)=-P_e(\nabla_uu+\nabla u^T\!\!\cdot
u)=-P_e(\nabla_uu)$ since $\nabla u^T\!\!\cdot
u=\operatorname{grad} (g(u,u))$.

It remains to show that $X_h(u)\in T^W_u\mathfrak{X}^s_{div}(M)$,
that is, there exists $d\in C^1_W(I,\mathfrak{X}^s_{div}(M))$ such
that $d(0)=u$ and $\dot{d}(0)=P_e(\nabla_uu)$. It suffices to
consider the curve $d(t):=P_e(T\eta(t)^T\circ u\circ
\eta(t))\in\mathfrak{X}^s_{div}(M)$  where $\eta(t)$ is a smooth
curve in $\mathcal{D}^s_\mu(M)$ such that $\eta(0)=id$ and
$\dot{\eta}(0)=u$. The fact that $d\in
C^1_W(I,\mathfrak{X}^s_{div}(M))$ is a consequence of that fact
that for $s>\frac{\operatorname{dim}(M)}{2}+1$ the map
\[
\omega_u : \mathcal{D}^s_\mu(M)\longrightarrow
\mathfrak{X}^s(M),\,\,\omega_u(\eta)=u\circ \eta
\]
is continuous and is $C^1$ as a map with values in
$\mathfrak{X}^{s-1}(M)$ (see the proof of Proposition 3.4 in
\cite{Eb1968}).$\;\;\;\;\blacksquare$

\medskip

Note that for $u\in C^1_W(I,\mathfrak{X}^s_{div}(M)),
s>\frac{\operatorname{dim}(M)}{2}+1$, and $h(u)=\frac{1}{2}\langle
u,u\rangle$ we have the following equivalent formulations of the
Euler equations
\begin{enumerate}
\item[{\rm (i)}]
\[
\partial_tu(t)+\nabla_{u(t)}u(t)=-\operatorname{grad} p(t)
\]
for some scalar function $p(t) : M\longrightarrow \mathbb{R}$,
\item[{\rm (ii)}]
\[
\partial_tu(t)=X_h(u(t)),
\]
\item[{\rm (iii)}] for all $f\in C^1_W(\mathfrak{X}^s_{div}(M))$
that admit a functional derivative $\delta f :
\mathfrak{X}^s_{div}(M)\longrightarrow \mathfrak{X}^s_{div}(M)$ we
have
\[
\left.\frac{d}{dt}\right|_{t=0}f(u(t))=\{f,h\}_+(u(t)).
\]
\end{enumerate}

It is well-known, in the standard case, that the flow of invariant
vector fields leaves the isotropy type submanifold invariant.
Using existence and uniqueness of the solution of the Euler
equations, we will show below that this property is still true for
the vector field $X_h : \mathfrak{X}^s_{div}(M)\longrightarrow
T^W\mathfrak{X}^s_{div}(M)$.

\begin{theorem}\label{conservation_of_isotropy} (Law of conservation of the
isotropy) Let $X_h$ be the Hamiltonian vector field associated
to the Euler equation and let $\widetilde{F}_t$ be its flow. Then:
\begin{enumerate}
\item[{\rm (i)}] $X_h\circ l_i=l_i\circ X_h$ and
$\widetilde{F}_t\circ l_i=l_i\circ \widetilde{F}_t$ for all $i\in
Iso^+$,

\item[{\rm (ii)}] for all closed subgroups $H$ of $Iso^+$,
$X_h^H:=X_h|_{\mathfrak{X}^s_{div}(M)_H}$ is a vector field on
$\mathfrak{X}^s_{div}(M)_H$, that is
\[
X_h^H : \mathfrak{X}^s_{div}(M)_H \longrightarrow
T^W\mathfrak{X}^s_{div}(M)_H,
\]
and the flow of $X^H_h$ is
\[
\widetilde{F}_t^H:=\widetilde{F}_t|_{\mathfrak{X}^s_{div}(M)_H} :
\mathfrak{X}^s_{div}(M)_H \longrightarrow
\mathfrak{X}^s_{div}(M)_H.
\]
\end{enumerate}
\end{theorem}
\textbf{Proof.} (i) The invariance of $X_h$ is a direct
computation using Lemma \ref{Iso^+invariance}. Let
$u_0\in\mathfrak{X}^s_{div}(M)$. We know that
$\widetilde{F}_t(l_i(u_0))$ is an integral curve of $X_h$ through
$l_i(u_0)$. One can check, using the invariance of $X_h$, that the
same is true for $l_i(\widetilde{F}_t(u_0))$. By uniqueness of the
solutions of Euler equations we obtain that $\widetilde{F}_t\circ
l_i=l_i\circ \widetilde{F}_t$, so we conclude that
$\widetilde{F}_t(\mathfrak{X}^s_{div}(M)_H)\subset\mathfrak{X}^s_{div}
(M)^H$. Using the bijectivity of the flow gives
$\widetilde{F}_t(\mathfrak{X}^s_{div}(M)_H)\subset\mathfrak{X}^s_{div}
(M)_H.$ (ii) is a direct verification.$\;\;\;\;\blacksquare$

\medskip

Note that the preceding theorem remains valid for other
$Iso^+$-invariant Hamiltonians $h$ such that their corresponding
motion equations
\[
\partial_tu(t)+\nabla_{\delta h(u(t))}u(t)+\nabla\delta h(u(t))^T\!\!\cdot u(t)
=-\operatorname{grad} p(t)
\]
admit a unique integral curve $u\in
C^1_W(I,\mathfrak{X}^s_{div}(M))$ for each initial condition.

\subsection{Poisson brackets on $\mathfrak{X}^s_{div}(M)_H$ and
$\mathfrak{X}^s_{div}(M)_H/N(H)$}

In this subsection we will define a Poisson bracket
$\{\,,\}_{\mathfrak{X}^s_H}$ on $\mathfrak{X}^s_{div}(M)_H$, and
we shall show in which sense the inclusion map $i_H :
\mathfrak{X}^s_{div}(M)_H\longrightarrow \mathfrak{X}^s_{div}(M)$
is a Poisson map. Then we will define a Poisson bracket
$\{\,,\}_{\mathfrak{X}^s_H/N(H)}$ on
$\mathfrak{X}^s_{div}(M)_H/N(H)$, and we shall show in which sense
the projection map $\pi_H : \mathfrak{X}^s_{div}(M)_H
\longrightarrow \mathfrak{X}^s_{div}(M)_H/N(H)$ is a Poisson map.

\begin{definition} For $k\geq 1, s>\frac{\operatorname{dim}(M)}{2}+1$,
and
$r\geq 0$ we define
\[
C^k_r(\mathfrak{X}^s_{div}(M)_H):=\left\{f\in
C^k(\mathfrak{X}^s_{div}(M)_H)\,|\,\exists\,\delta f :
\mathfrak{X}^s_{div}(M)_H\longrightarrow
\mathfrak{X}^r_{div}(M)^H\right\}.
\]
The \textbf{Poisson bracket} of $f,g\in
C^k_r(\mathfrak{X}^s_{div}(M)_H),
r>\frac{\operatorname{dim}(M)}{2}+1$, is defined by
\[
\{f,g\}_{\mathfrak{X}^s_H}(u):=\langle u,[\delta g(u),\delta
f(u)]\rangle.
\]
\end{definition}

\begin{theorem}\label{i_H is Poisson} Let $k\geq 1$ and
$r,s>\frac{\operatorname{dim}(M)}{2}+1$.
\begin{enumerate}
\item[{\rm (i)}]For all $f\in C^k_r(\mathfrak{X}^s_{div}(M))^H$,
we have
\[
f\circ i_H\in C^k_r(\mathfrak{X}^s_{div}(M)_H)\;\text{ and }\;\
\delta(f\circ i_H)=\delta f\circ i_H.
\]
\item[{\rm (ii)}]For all $f,g\in
C^k_r(\mathfrak{X}^s_{div}(M))^H$, we have
\[
\{f\circ i_H,g\circ i_H\}_{\mathfrak{X}^s_H}=\{f,g\}_+\circ
i_H\;\text{ on }\;\mathfrak{X}^s_{div}(M)_H.
\]
\end{enumerate}
\end{theorem}
\textbf{Proof.} (i) For all $u\in\mathfrak{X}^s_{div}(M)_H$ and
for all $v\in\mathfrak{X}^s_{div}(M)^H$, we have:
\[
D(f\circ
i_H)(u)(v)=Df(i_H(u))(Di_H(u)(v))=Df(i_H(u))(i_H(v))=\langle
\delta f(i_H(u)),i_H(v)\rangle.
\]
We now show that $\delta f(i_H(u))\in\mathfrak{X}^r_{div}(M)^H$.
If $i\in H$, we have
\begin{align*}
l_i(\delta f(u))&=l_i(\delta f(l_{i^{-1}}(u))),\;\text{ since
}\;u\in\mathfrak{X}^s_{div}(M)_H\\
&=\delta(f\circ l_{i^{-1}})(u),\;\text{ by Lemma
\ref{Iso^+invariance}}\\
&=\delta f(u),\;\text{ since }\;f\in
C^k_r(\mathfrak{X}^s_{div}(M))^H.
\end{align*}
Thus we obtain the existence of $\delta(f\circ i_H)=\delta f\circ
i_H$.\\
(ii) Since $f\circ i_H, g\circ i_H\in
C^k_r(\mathfrak{X}^s_{div}(M)_H)$ we can compute, for
$u\in\mathfrak{X}^s_{div}(M)_H$,
\begin{align*}
\{f\circ i_H, g\circ i_H\}_{\mathfrak{X}^s_H}(u)&=\langle
u,[\delta (g\circ
i_H)(u),\delta (f\circ i_H)(u)]\rangle\\
&=\langle u,[\delta g(i_H(u)),\delta f(i_H(u))]\rangle\;\;\text{ by (i)}\\
&=\langle i_H(u),[\delta g(i_H(u)),\delta f(i_H(u))]\rangle\\
&=\{f,g\}_+(i_H(u)).\;\;\;\;\blacksquare
\end{align*}

In order to define the Poisson bracket on
$\mathfrak{X}^s_{div}(M)_H/N(H)$, we will need the following
subsets of $C^1_W(\mathfrak{X}^s_{div}(M)_H)$ and
$C^1_W(\mathfrak{X}^s_{div}(M)_H/N(H))$.

\begin{definition}\label{functional_derivative}
\begin{enumerate} \item[{\rm
(i)}] For $r\geq 0,s\geq 1$ we define the set
\[
C^1_{Wr}(\mathfrak{X}^s_{div}(M)):=\left\{f\in
C^1_W(\mathfrak{X}^s_{div}(M))\,|\,\exists\,\delta f :
\mathfrak{X}^s_{div}(M)\longrightarrow
\mathfrak{X}^r_{div}(M)\right\}
\]
where $\delta f$ is the functional derivative of $f$ with respect
to $\langle\,,\rangle$, that is,
\[
\langle \delta f(u),v\rangle =Df(u)(v),\,\,\text{ for all }
u,v\in\mathfrak{X}^s_{div}(M).
\]
This is possible since $C^1_W(\mathfrak{X}^s_{div}(M))\subset
C^1(\mathfrak{X}^s_{div}(M))$.

\item[{\rm (ii)}] In a similar way we define
\[
C^1_{Wr}(\mathfrak{X}^s_{div}(M)_H):=\left\{f\in
C^1_W(\mathfrak{X}^s_{div}(M)_H)\,|\,\exists\,\delta f :
\mathfrak{X}^s_{div}(M)_H\longrightarrow
\mathfrak{X}^r_{div}(M)^H\right\}.
\]

\item[{\rm (iii)}] For $r\geq 0,
s>\frac{\operatorname{dim}(M)}{2}+2$ we define the set
\[
C^1_{Wr}(\mathfrak{X}^s_{div}(M)_H/N(H)):=\left\{\varphi\in
C^1_W(\mathfrak{X}^s_{div}(M)_H/N(H))\,|\,\varphi\circ\pi_H\in
C^1_{Wr}(\mathfrak{X}^s_{div}(M)_H)\right\},
\]
where $\pi_H :\mathfrak{X}^s_{div}(M)_H\longrightarrow
\mathfrak{X}^s_{div}(M)_H/N(H)$.
\end{enumerate}
\end{definition}

Recall that on $\mathfrak{X}^s_{div}(M)_H/N(H),
s>\frac{\operatorname{dim}(M)}{2}+2$, we can define a Riemannian
metric
\[
\gamma^H(\omega)(\xi_\omega,\eta_\omega):
=\langle\operatorname{Hor}^H_{u_\omega}(\xi_\omega),
\operatorname{Hor}^H_{u_\omega}(\eta_\omega)\rangle
\]
where $\omega\in\mathfrak{X}^s_{div}(M)_H/N(H)$, $
\xi_\omega,\eta_\omega\in
T^W_\omega(\mathfrak{X}^s_{div}(M)_H/N(H))$, and $u_\omega$ is any
element in $\pi^{-1}_H(\omega)$.

So for $\varphi\in C^1_{Wr}(\mathfrak{X}^s_{div}(M)_H/N(H))$ we
have, using Theorem \ref{Horizontalchainrule},
\begin{align*}
d^W\varphi(\omega)(\xi_\omega)&
=D^W(\varphi\circ\pi_H)(u_\omega)(\operatorname{Hor}^H_{u_\omega}(\xi_\omega))
\\
&=\langle\delta
(\varphi\circ\pi_H)(u_\omega),\operatorname{Hor}^H_{u_\omega}(\xi_\omega)
\rangle\\
&=\langle\operatorname{Hor}^H_{u_\omega}(T^W_{u_\omega}\pi_H(\delta
(\varphi\circ\pi_H)(u_\omega))),\operatorname{Hor}^H_{u_\omega}(\xi_\omega)
\rangle\\
&=\gamma^H(\omega)(T^W_{u_\omega}\pi_H(\delta
(\varphi\circ\pi_H)(u_\omega)),\xi_\omega).
\end{align*}
We conclude that any $\varphi\in
C^1_{Wr}(\mathfrak{X}^s_{div}(M)_H/N(H))$ admits a functional
derivative with respect to $\gamma^H$. It is given by
\[
\delta \varphi :\mathfrak{X}^s_{div}(M)_H/N(H)\longrightarrow
T^W(\mathfrak{X}^s_{div}(M)_H/N(H)),\,\,\delta\varphi(\omega)=T^W_
{u_\omega}\pi_H(\delta(\varphi\circ\pi_H)(u_\omega)),
\]
where $u_\omega$ is any element in $\pi_H^{-1}(\omega)$. Since
$\varphi\circ\pi_H$ is $N(H)$-invariant,
$\delta(\varphi\circ\pi_H)(\omega)$ is horizontal by \eqref{weak
derivative horizontal}, so we have
\[
\delta(\varphi\circ\pi_H)(u_\omega)=\operatorname{Hor}^H_{u_\omega}
(\delta\varphi
(\omega)).
\]

\begin{lemma}\label{LieBracket} For all
$u,v\in\mathfrak{X}^s_{div}(M),
s>\frac{\operatorname{dim}(M)}{2}+1$, we have
\[
[u,v]\in\mathfrak{X}^s_{div}(M)^W
\]
\end{lemma}
\textbf{Proof.} It suffices to consider the curve
$d(t):=\eta(t)^*v$ where $\eta(t)$ is a smooth curve in
$\mathcal{D}^s_\mu(M)$ such that $\eta(0)=id$ and
$\dot{\eta}(0)=u$. So we have $d\in
C^1_W(I,\mathfrak{X}^s_{div}(M))$ and
$\dot{d}(0)=L_uv=[u,v].\;\;\;\;\blacksquare$

\begin{definition}
\label{def poisson} The \textbf{Poisson bracket} of
$\varphi,\psi\in C^1_{Wr} (\mathfrak{X}^s_{div}(M)_H/N(H)),$
$r,s>\frac{\operatorname{dim}(M)}{2}+2$, is defined by
\[
\{\varphi,\psi\}_{\mathfrak{X}^s_H/N(H)}(\omega):=\gamma^H(\omega)\left(S^H
(\omega),[\![\delta\psi
(\omega),\delta\varphi(\omega)]\!]^H\right),
\]
where\begin{enumerate}

\item[\rm (i)]$S^H : \mathfrak{X}^s_{div}(M)_H/N(H)\longrightarrow
T^W(\mathfrak{X}^s_{div}(M)_H/N(H))$ is the vector field on
$\mathfrak{X}^s_{div}(M)_H/N(H)$ defined by
$S^H(\omega):=T^W_{u_\omega}\pi_H(u_\omega)$ for any
$u_\omega\in\pi_H^{-1}(\omega)$ and

\item[\rm (ii)] $[\![\,,]\!]^H$ is the \textbf{reduced Lie
bracket} defined by
\[
[\![\xi_\omega,\eta_\omega]\!]^H:=T^W_{u_\omega}\pi_H\left([\operatorname
{Hor}^H_{u_\omega}(\xi_\omega),\operatorname{Hor}^H_{u_\omega}(\eta_\omega)]
\right).
\]
\end{enumerate}
\end{definition}

Note that since
$\operatorname{Hor}^H_{u_\omega}(\delta\varphi(\omega))
=\delta(\varphi\circ\pi_H)(u_\omega)\in\mathfrak{X}^r_{div}(M),
r>\frac{\operatorname{dim}(M)}{2}+2$, we have
\[
[\operatorname{Hor}^H_{u_\omega}(\delta\varphi(\omega)),\operatorname{Hor}^H_
{u_\omega}(\delta\psi(\omega))]\in\mathfrak{X}^r_{div}(M)^W
\]
by Lemma \ref{LieBracket}. So
$T^W_{u_\omega}\pi_H([\operatorname{Hor}_{u_\omega}(\delta\varphi
(\omega)),\operatorname{Hor}_{u_\omega}(\delta\psi(\omega))])$ is
well-defined.

\begin{theorem} \label{pi_is_Poisson} For
$s>\frac{\operatorname{dim}(M)}{2}+2
$, the projection
\[
\pi_H :\mathfrak{X}^s_{div}(M)_H\longrightarrow
\mathfrak{X}^s_{div}(M)_H/N(H)
\]
is a Poisson map, that is, for all $\varphi,\psi\in
C^1_{Wr}(\mathfrak{X}^s_{div}(M)_H/N(H)),
r>\frac{\operatorname{dim}(M)}{2}+2$, we have
\[
\{\varphi\circ\pi_H,\psi\circ\pi_H\}_{\mathfrak{X}^s_H}=\{\varphi,\psi\}_
{\mathfrak{X}^s_H/N(H)}\circ\pi_H.
\]
\end{theorem}
\textbf{Proof.} Since $\varphi,\psi\in
C^1_{Wr}(\mathfrak{X}^s_{div}(M)_H/N(H))$, we have
$\varphi\circ\pi_H,\psi\circ\pi_H\in
C^1_{Wr}(\mathfrak{X}^s_{div}(M)_H)\subset
C^1_r(\mathfrak{X}^s_{div}(M)_H)$ by Definition
\ref{functional_derivative} (iii). Note that $u\in
T^W_u\mathfrak{X}^s_{div}(M)_H$ is horizontal, since we have
$\langle u,\xi_{\mathfrak{X}^s_{div}(M)}(u) \rangle = \langle
u,[u,X_\xi]\rangle=\langle u,\nabla_uX_\xi\rangle-\langle
u,\nabla_{X_\xi}u\rangle=0$. Indeed, integrating by parts we have
$\langle u,\nabla_{X_\xi}u\rangle=-\langle
\nabla_{X_\xi}u,u\rangle$ so $\langle u,\nabla_{X_\xi}u\rangle=0$.
Since $X_\xi$ is a Killing vector field, $\nabla X_\xi$ is skew
symmetric, so $\langle u,\nabla_uX_\xi\rangle=0$. This implies
that $u = \operatorname{Hor}^H_u (S^H(\pi_H(u))) $. So
for $u\in\mathfrak{X}^s_{div}(M)_H$ we have
\begin{align*}
\{\varphi\circ\pi_H,&\psi\circ\pi_H\}_{\mathfrak{X}^s_H}(u)=\langle
u,[\delta(\psi\circ\pi_H)(u),\delta(\varphi\circ\pi_H)(u)]\rangle\\
&=\langle
\operatorname{Hor}^H_u(S^H(\pi_H(u))),[\delta(\psi\circ\pi_H)(u),\delta
(\varphi\circ\pi_H)(u)]\rangle\\
&=\langle
\operatorname{Hor}^H_u(S^H(\pi_H(u))),\operatorname{hor}^H_u([\delta
(\psi\circ\pi_H)(u),\delta(\varphi\circ\pi_H)(u)])\rangle\\
&=\langle
\operatorname{Hor}^H_u(S^H(\pi_H(u))),\operatorname{Hor}^H_u(T^W_u\pi_H([\delta
(\psi\circ\pi_H)(u),\delta(\varphi\circ\pi_H)(u)]))\rangle\\
&=\gamma^H(\pi_H(u))\left(S^H(\pi_H
(u)),T^W_u\pi_H([\operatorname{Hor}^H_u(\delta\psi
(\pi_H(u))),\operatorname{Hor}^H_u(\delta\varphi(\pi_H(u)))])\right)\\
&=\{\varphi,\psi\}_{\mathfrak{X}^s_H/N(H)}(\pi_H(u))
\end{align*}
by Definition \ref{def poisson}.\;\;\;\;$\blacksquare$

\subsection{The reduced Euler equations}

Recall that the Hamiltonian of the Euler equations is given by
$h(u)=\frac{1}{2}\langle u,u\rangle$. We have $h\in
C^1_{Ws}(\mathfrak{X}^s_{div}(M))$. The Hamiltonian vector field
of $h$ with respect to the Poisson bracket $\{\,,\}_+$, is given
by (see \eqref{Poissonformulation})
\[
X_h : \mathfrak{X}^s_{div}(M)\longrightarrow
T^W\mathfrak{X}^s_{div}(M),\;\; X_h(u)=(u,-P_e(\nabla_uu)).
\]
By restriction to the isotropy type manifold, we have $h\in
C^1_{Ws}(\mathfrak{X}^s_{div}(M)_H)$ and one can show that the
Hamiltonian vector field with respect to the Poisson bracket
$\{\,,\}_{\mathfrak{X}^s_H}$ is given by the restriction
$X^H_h:=X_h|\mathfrak{X}^s_{div}(M)_H$, that is,
\[
X_h^H : \mathfrak{X}^s_{div}(M)_H\longrightarrow
T^W\mathfrak{X}^s_{div}(M)_H,\;\; X_h(u)=(u,-P_e(\nabla_uu)).
\]
By $N(H)$-invariance, $h$ induces a unique function
$\mathfrak{h}^H :
\mathfrak{X}^s_{div}(M)_H/N(H)\longrightarrow\mathbb{R}$, called
the reduced Hamiltonian, such that $h=\mathfrak{h}^H\circ\pi_H$ on
$\mathfrak{X}^s_{div}(M)_H$. By definition, we have
$\mathfrak{h}^H\in C^1_{Ws}(\mathfrak{X}^s_{div}(M)_H/N(H))$ and a
direct computation gives
\[
\mathfrak{h}^H(\omega)=\frac{1}{2}\gamma^H(\omega)(S^H(\omega),S^H(\omega)).
\]
Likewise, the Hamiltonian vector field $X_h^H$ induces a unique
vector field
\[
X_{\mathfrak{h}^H} : \mathfrak{X}^s_{div}(M)_H/N(H)\longrightarrow
T^W(\mathfrak{X}^s_{div}(M)_H/N(H))
\]
such that $X_{\mathfrak{h}^H}\circ\pi_H=T^W\pi_H\circ X_h^H$ on
$\mathfrak{X}^s_{div}(M)_H$.  More precisely, we have the following
result.

\begin{lemma} For $s>\frac{\operatorname{dim}(M)}{2}+2$, $X_{\mathfrak{h}^H}
:
\mathfrak{X}^s_{div}(M)_H/N(H)\longrightarrow
T^W(\mathfrak{X}^s_{div}(M)_H/N(H))$ is given by
\[
X_{\mathfrak{h}^H}(\omega)=-\mathcal{P}^H_e(\nabla^H_{S^H(\omega)}S^H(\omega)),
\]
where
\begin{enumerate}

\item[\rm (i)] $\mathcal{P}^H_e :
T^W(\mathfrak{X}^s(M)_H/N(H))\longrightarrow
T^W(\mathfrak{X}^s_{div}(M)_H/N(H))$ is the \textbf{reduced Hodge
projector} given by
\[
\mathcal{P}^H_e(\xi_\omega):=T^W_{u_\omega}\pi_H(P_e(\operatorname{Hor}^H_
{u_\omega}
(\xi_\omega)))
\]
for any $u_\omega \in \pi_H^{-1}(\omega)$ and

\item[\rm (ii)] for $\xi_\omega, \eta_\omega\in
T^W_\omega(\mathfrak{X}^s_{div}(M)_H/N(H))$,
\[
\nabla^H_{\eta_\omega}\xi_\omega:=T^W_{u_\omega}\pi_H(\nabla_{\operatorname
{Hor}^H_{u_\omega}(\eta_\omega)}\operatorname{Hor}^H_{u_\omega}
(\xi_\omega))
\]
for any $u_\omega \in \pi_H^{-1}(\omega)$.
\end{enumerate}
\end{lemma}
\textbf{Proof.} For any $u\in\pi^{-1}_H(\omega)$ we have
\begin{align*}X_{\mathfrak{h}^H}(\omega)&=T^W_u\pi_H(X_h(u))\\
&=-T^W_u\pi_H(P_e(\nabla_uu))\\
&=-T^W_u\pi_H(P_e(\operatorname{hor}^H_u(\nabla_uu)+\operatorname{ver}^H_u
(\nabla_uu)))\\
&=-T^W_u\pi_H(P_e(\operatorname{hor}^H_u(\nabla_uu))+\operatorname{ver}^H_u
(\nabla_uu))\\
&=-T^W_u\pi_H(P_e(\operatorname{Hor}^H_u(T_u^W\pi_H(\nabla_uu))))\\
&=-\mathcal{P}_e^H(T^W_u\pi_H(\nabla_uu))\\
&=-\mathcal{P}_e^H(T^W_u\pi_H(\nabla_{\operatorname{Hor}^H_u(S(\omega))}
\operatorname{Hor}^H_u(S^H(\omega))))\\
&=-\mathcal{P}_e^H(\nabla^H_{S^H(\omega)}S^H(\omega)).
\end{align*}
For the fourth equality we use that $P_e$ is the identity on
vertical vector fields, since they belong to
$\mathfrak{X}^s_{div}(M).\;\;\;\;\blacksquare$

\medskip

One can show as above that $X_{\mathfrak{h}^H}$ is the Hamiltonian
vector field associated to $\mathfrak{h}^H$ with respect to
$\{\,,\}_{\mathfrak{X}^s_H/N(H)}$, that is,
\[
d^W\varphi(\omega)(X_{\mathfrak{h}^H}(\omega))=\{\varphi,\mathfrak{h}^H\}_
{\mathfrak{X}^s_H/N(H)}(
\omega),\,\,\text{ for all }\varphi\in
C^1_{Wr}(\mathfrak{X}^s_{div}(M)_H/N(H)).
\]

Let $u(t)$ be an integral curve of the Euler equations. Thus we
have $u\in C^1_W(I,\mathfrak{X}^s_{div}(M))$. If the initial
condition $u(0)=u_0$ is in $\mathfrak{X}^s_{div}(M)_H$ then by the
law of conservation of isotropy (Theorem
\ref{conservation_of_isotropy}) we have $u\in
C^1_W(I,\mathfrak{X}^s_{div}(M)_H)$. Let $\omega:=\pi_H\circ u$;
then $\omega \in C^1_W(I,\mathfrak{X}^s_{div}(M)_H/N(H))$ and
$\left.\frac{d}{dt}\right|_{t=0}\omega(t)=X_{\mathfrak{h}^H}(\omega(t))$,
that is,
\[
\left.\frac{d}{dt}\right|_{t=0}\omega(t)=-\mathcal{P}^H_e(\nabla^H_{S^H(\omega
(t))}S^H(\omega(t))).
\]
These equations are called the \textbf{reduced Euler equations} on
$\mathfrak{X}^s_{div}(M)_H/N(H)$.

\medskip

Let $\widetilde{F}_t$ be the flow of the Euler equations on
$\mathfrak{X}^s_{div}(M)$ and  $\widetilde{F}_t^H$  the flow of
the Euler equations on $\mathfrak{X}^s_{div}(M)_H$. Define
$\widetilde{\widetilde{F}}_t\,^H(\omega):=\pi_H\left(\widetilde{F}_t^H
(u_\omega)\right)$
where $u_\omega$ is such that $\pi_H(u_\omega)=\omega$. Then
$\widetilde{\widetilde{F}}_t\,^H(\omega)$ is the flow of the
reduced Euler equations on $\mathfrak{X}^s_{div}(M)_H/N(H)$, that
is, $\omega(t):=\widetilde{\widetilde{F}}_t\,^H(\omega_0)$ is the
integral curve through $\omega_0$. Moreover, we have the
commutative diagram
$$\xymatrix{
\mathfrak{X}^s_{div}(M)_H \ar[r]^{\widetilde{F}_t^H}
\ar[d]_{\pi_H}&
\mathfrak{X}^s_{div}(M)_H \ar [d]^{\pi_H}\\
\mathfrak{X}^s_{div}(M)_H/N(H)\ar[r]^{\widetilde{\widetilde{F}}_t\,^H}
&
\mathfrak{X}^s_{div}(M)_H/N(H).\\
}$$

We already know that $\widetilde{F}_t^H$ and $\pi_H$ are Poisson
maps in the precise sense given in Theorems \ref{VaMa} and
\ref{pi_is_Poisson}. In the following theorem we show in which sense
$\widetilde{\widetilde{F}}_t\,^H$ is a Poisson map.

\begin{theorem}
\label{flow_is_Poisson} For all $\varphi,\psi\in
C^1_{Wr}(\mathfrak{X}^s_{div}(M)_H/N(H)),
r,s>\frac{\operatorname{dim}(M)}{2}+2$, such that $s+1\geq r$ we
have
\[
\{\varphi,\psi\}_{\mathfrak{X}^s_H/N(H)}\circ\widetilde{\widetilde{F}}_t\,^H
=\{\varphi\circ\widetilde{\widetilde{F}}_t\,^H,\psi\circ\widetilde{\widetilde
{F}}
_t\,^H\}_{\mathfrak{X}^s_H/N(H)}
\]
on $\mathfrak{X}^{s+2}_{div}(M)_H/N(H)$.
\end{theorem}
\textbf{Proof.} For $\omega\in\mathfrak{X}^{s+2}_{div}(M)_H/N(H)$
we have for any $u_\omega\in\pi_H^{-1}(\omega)$
\begin{align*}
\{\varphi,\psi\}_{\mathfrak{X}^s_H/N(H)}&(\widetilde{\widetilde
{F}}_t\,^H(\omega))
=\{\varphi,\psi\}_{\mathfrak{X}^s_H/N(H)}(\pi_H(\widetilde{F}_t^H(u_\omega)))
\\
&=\{\varphi\circ\pi_H,\psi\circ\pi_H\}_{\mathfrak{X}^s_H}(\widetilde{F}_t^H
(u_\omega))\,\,\text{
by Theorem \ref{pi_is_Poisson}}\\
&=\{\varphi\circ\pi_H\circ\widetilde{F}_t^H,\psi\circ\pi_H\circ\widetilde{F}
_t^H\}_{\mathfrak{X}^s_H}
(u_\omega)\\
&=\{\varphi\circ\widetilde{\widetilde{F}}_t\,^H\circ\pi_H,\psi\circ\widetilde
{\widetilde{F}}_t\,^H\circ\pi_H\}_{\mathfrak{X}^s_H}(u_\omega)\,\,\text{
by the commutative diagram}\\
&=\{\varphi\circ\widetilde{\widetilde{F}}_t\,^H,\psi\circ\widetilde{\widetilde
{F}}
_t\,^H\}_{\mathfrak{X}^s_H/N(H)}(\pi_H(u_\omega))\,\,\text{
by Theorem \ref{pi_is_Poisson}}\\
&=\{\varphi\circ\widetilde{\widetilde{F}}_t\,^H,\psi\circ\widetilde{\widetilde
{F}}
_t\,^H\}_{\mathfrak{X}^s_H/N(H)}(\omega).
\end{align*}

For the third equality we use Theorem \ref{VaMa} (iii) which is
applicable since we assumed that  $s+1\geq r$ and
$u_\omega\in\mathfrak{X}^{s+2}_{div} (M)_H$.

For the fifth equality we use that
$\varphi\circ\widetilde{\widetilde{F}}_t\,^H\in
C^1_{Wr}(\mathfrak{X}^{s+2}_{div}(M)_H)$. Namely, from
\cite{VaMa2005} we know that for all $f\in
C^1_r(\mathfrak{X}^s_{div}(M))$ we have $f\circ\widetilde{F}_t\in
C^1_r(\mathfrak{X}^{s+1}_{div}(M))$. So we obtain that
$f\circ\widetilde{F}_t^H|_{\mathfrak{X}^{s+2}_{div}(M)_H}\in
C^1_{Wr}(\mathfrak{X}^{s+2}_{div}(M)_H)$ and
$\varphi\circ\widetilde{\widetilde{F}}
_t\,^H\circ\pi_H=\varphi\circ\pi_H\circ\widetilde
{F}_t^H\in C^1_{Wr}(\mathfrak{X}^{s+2}_{div}(M)_H)$ because
$\varphi\circ\pi_H\in C^1_{Wr}(\mathfrak{X}^s_{div}(M)_H)\subset
C^1_r(\mathfrak{X}^s_{div}(M)_H)$. Thus, by definition, we find that
$\varphi\circ\widetilde{\widetilde{F}}_t\,^H\in
C^1_{Wr}(\mathfrak{X}^{s+2}_{div}(M)_H).\;\;\;\;\blacksquare$
\medskip

Finally we obtain the following commuting diagram in which all
maps are Poisson in the precise sense given in Theorems
\ref{VaMa}, \ref{i_H is Poisson}, \ref{pi_is_Poisson}, and
\ref{flow_is_Poisson}.

$$\xymatrix{
T\mathcal{D}^s_\mu(M) \ar[r]^{F_t} \ar[d]_{\pi_R}&
T\mathcal{D}^s_\mu(M)\ar [d]^{\pi_R}\\
\mathfrak{X}^s_{div}(M) \ar[r]^{\widetilde{F}_t} &
\mathfrak{X}^s_{div}(M) \\
\mathfrak{X}^s_{div}(M)_H \ar[u]^{i_H} \ar[r]^{\widetilde{F}_t^H}
\ar[d]_{\pi_H}&
\mathfrak{X}^s_{div}(M)_H \ar[u]^{i_H}\ar [d]^{\pi_H}\\
\mathfrak{X}^s_{div}(M)_H/N(H)\ar[r]^{\widetilde{\widetilde{F}}_t\,^H}
&
\mathfrak{X}^s_{div}(M)_H/N(H).\\
}$$

The same results can be obtained for the averaged Euler equations
by using \cite{GBRa2005} instead of \cite{VaMa2005}.
\medskip

\noindent {\bf Acknowledgments.} We thank J. Marsden for
suggesting the problem and for several discussions during the
elaboration of the paper. Our thanks to P. Michor for
conversations regarding infinite dimensional geometry. The first
author was fully supported by a doctoral fellowship of the EPFL.
The second author acknowledges the partial support of the  Swiss
NSF.

{\footnotesize

\bibliographystyle{new}

}

\end{document}